\documentclass{amsart}
\usepackage{amssymb, amsmath, amsfonts, amscd}

\input xypic

\theoremstyle{plain}
\newtheorem{theorem}{Theorem}[section]
\newtheorem{corollary}[theorem]{Corollary}
\newtheorem{lemma}[theorem]{Lemma}
\newtheorem{proposition}[theorem]{Proposition}

\newtheorem{example}[theorem]{Example}

\theoremstyle{definition}
\newtheorem{definition}[theorem]{Definition}

\theoremstyle{remark}

\numberwithin{equation}{theorem}

\newcommand{\Hom}{\operatorname{Hom} }
\newcommand{\onabla}{\overline{\nabla} }
\newcommand{\End}{\operatorname{End} }
\newcommand{\Ext}{\operatorname{Ext} }

\renewcommand{\H}{\operatorname{H} }

\newcommand{\Spec}{\operatorname{Spec} }

\newcommand{\Der}{\operatorname{Der} }

\newcommand{\Conn}{\operatorname{Conn} }
\newcommand{\Diff}{\operatorname{Diff} }

\newcommand{\D}{\operatorname{D} }

\newcommand{\K}{\operatorname{K}}

\newcommand{\Z}{\operatorname{Z}} 
\newcommand{\ext}{\operatorname{ext}}

\begin{document}

\title{A universal  characteristic class for vector bundles with a connection}

\author{Helge Maakestad }



\date{November 2024}

\begin{abstract} The Whitney sum formula implies that the Chern class $c_i(E)$ for a stably trivial vector bundle $E$ on an affine scheme $X$ is trivial. The formula also implies that the class in the Grothendieck group  $[E] \in \operatorname{K}_0(X)$ is trivial. If $X$ is an affine regular scheme of finite type over a field $k$ and if $E$ is the tangent or cotangent bundle, it follows $E$ is stably trivial, hence $c_i(E)=0$ and $[E]=rk(E)$. Hence there is a large class of non trivial vector bundles with trivial Chern classes and trivial classes in the Chow group and Grothendieck group. The aim of this paper is to intruduce a new characteristic class $c(E)$ that is non trivial for stably trivial vector bundles. The class $c(E)$ lives in a pointed cohomology torsor and most classical characteristic classes live in a cohomology group, hence the class $c(E)$ is a new characteristic class not previously studied. The class $c(E)$ is universal in the sense that it can be used to construct and study a large class of characteristic classes. In this paper i prove the existence of the class $c(E)$ and prove that it is non trivial for the tangent bundle on the real two sphere. The tangent bundle on the real two sphere is stably trivial
hence its Chern classes are trivial. Hence the class $c(E)$ is stronger than the Chern classes. The class may be interesting in the study of the "cancellation problem" in affine algebraic geometry.

\end{abstract}

\maketitle
\tableofcontents

\section{Introduction}   The aim of this paper is to introduce new methods in the study of connections and higher order differential operators on modules and to use these methods to introduce a new characteristic class $c(W)$ - the fundamental class of $W$ - for a finite rank vector bundle $W$  on an affine scheme $X$. The characteristic class is not a characteristic class in the "classical sense" since the class lives in a cohomology torsor.  Most characteristic classes live in a cohomology group. It is universal in the sense that any characteristic class $u(W)$ that can be defined in terms of the curvature of a connection can be constructed from the class $c(W)$. Hence the class $c(W)$ can be used to study any such theory of characteristic classes. In fact there is a classical construction of Teleman, which construct for any  relative tangent sequence $ 0 \rightarrow W \rightarrow  L_1 \rightarrow L \rightarrow 0$
associated to any principal fiber bundle $\pi: P_1 \rightarrow S$ and any invariant polynomial $g$ and any connection $\gamma$, a characteristic class
$[g(\Omega_{\gamma})]$ which is independent of choice of connection $\gamma$. If $c(W)$ is the fundamental class of $W$  and if $c(W)$ is trivial it follows $[g(\Omega_{\gamma})]$ is trivial.
Hence the fundamental class $c(W)$ can be used to construct and study the class $[g(\Omega_{\gamma})]$ of Teleman. Kong constructs in \cite{kong} the Euler class $e(P,h)$ of an oriented projective module $(P,h)$ using the Pfaffian
of the curvature of a connection on $P$ compatible with the orientation $h$. The class $e(P,h)$ is independent of choice of connection and is trivial of $P$ has a flat connection $\nabla$ compatible with $h$. Hence the class $e(P,h)$ can be constructed and studied using $c(P)$. 

The methods introduced in this paper enable us to calculate explicit non-trivial examples "by hand" and I do this in the case ot the cotangent bundle on the real two sphere in Theorem \ref{maincalc}.
I get an example of a "naturally ocurring" non flat connection that is not of curvature type $f$ where $f$ is a 2-cocycle. I also prove that the class $c(\Omega)$ is non trivial, where $\Omega$ is the cotangent bundle on the real two sphere.

The motivation for the introduction of this characteristic class is the following: There is  a large class of finite rank vector bundles on an affine scheme that are stably trivial. If $X$ is an affine regular complete intersection it follows the cotangent bundle $\Omega^1_X$ and tangent bundle $T_X$ are stably trivial. 
It follows $\Omega^1_X$ and $T_X$ have vanishing Chern classes.  If $E$ is a stably trivial finite rank vector bundle on $X$, it follows its class in the Grothendieck group $\operatorname{K}(X)$ is trivial: There is an equality of classes $[E]= d[\mathcal{O}_X]$ in $\K(X)$ where $d$ is an integer. Hence there is a large class of finite rank non trivial vector bundles that have vanishing Chern classes and trivial classes in the Chow group and Grothendieck group and there is need for a characteristic class that is non trivial for such bundles.   The search for a new characteristic class that is non trivial for such vector bundles was inpired by the following general result on explicit formulas for the curvature of a connection on a projective module. The result shows that if you choose a connection $\nabla$ on a vector bundle $E$ it follows $\nabla$ will in general be non-flat and that  it is a non trivial problem to determine if a vector bundle $E$ has a flat connection:

Let $A$ be any commutative unital ring and let  $E$ be a finite rank projective $A$-module. Let $p: A^n \rightarrow E, s: E \rightarrow A^n$ be $A$-linear maps with $p \circ s = Id_E$ and $n \geq 1$ an integer.  Let $\phi:=s \circ p$  be the idempotent element $\phi \in \End_A(A^n)$ defining the module $E$. Using the split surjection $(p,s)$ I define a connection $\nabla_E$ on $E$ and an $(L, \phi)$-connection $\onabla_E$ on $A^n$ with the following properties:

There is a map of $A$-modules $\pi_A: \End_A(A^n) \rightarrow \End_A(E)$ defined by $\pi_A( z):= p \circ z \circ s$  and  a commutative diagram

\[
\diagram          &     \End_A(A^n) \dto^{\pi_A}    \\
                      \wedge^2 L \urto^{R_{\onabla_E}} \rto^{R_{\nabla_E}}  &   \End_A(E)    
\enddiagram.
\]

\begin{theorem}\label{curvatureintro}  Let $R_{\nabla_E}$ be the curvature of $\nabla_E$ and let  $R_{\onabla_E}$ be the curvature of $\onabla_E$. For any elements $x,y \in \Der_k(A)$ it follows 

\[ R_{\onabla_E}(x,y) = \phi \circ [x(\phi), y(\phi)] .\]

We recover the curvature $R_{\nabla_E}(x,y)$ via the formula

\[  R_{\nabla_E} = \pi_A \circ  R_{\onabla_E}  .\]

\end{theorem}

Here $x(\phi)$ and $y(\phi)$ are endomorphisms of $A^n$ defined using a basis $B$ of $A^n$ and the idempotent element $\phi$. The product $[x(\phi), y(\phi)]$ is the Lie product in the associative ring $\End_A(A^n)$. 
Hence we have "lifted" the curvature $R_{\nabla_E}$ to the curvature of the $(L, \phi)$-connection $\onabla_E$ on $A^n$. If we choose a basis $B$ for $A^n$ and a matrix $[\phi]_B=[a_{ij}]$ of $\phi$ in the basis $B$, 
it follows there is an equality of matrices

\begin{align}
&\label{curv}[R_{\onabla_E}(x,y)]_B =  [ a_{ij}]( [x(a_{ij})] [y(a_{ij})] - [y(a_{ij}) ][x(a_{ij})] ).
\end{align}

This is proved in Theorem \ref{curvature} and Corollary \ref{maincorollary}. Hence if we choose an arbitrary finite rank projective $A$-module $E$ with a split surjection $(p,s)$ defining $E$, it follows the corresponding connection $\nabla_E$ will in general be non flat. It is a non trivial problem
to determine if  $E$ has a flat algebraic connection. Given an arbitrary connection $\nabla_1$ on $E$, there is always a potential $P$ and an equality $\nabla_1 = \nabla_E + P$. Hence the formula $\ref{curv}$ expresses the curvature of   $\nabla_1$ in terms of the idempotent element $\phi$ and the potential $P$. We get an explicit formula for the curvature of any connection $\nabla_1$. The formula implies a classical formula for the Chern character of a projective module with values in cyclic homology (see \cite{loday}, Theorem 8.3.2)

Kunz has studied similar questions for vector bundles 
on affine "Riemannian schemes". Using Kunz' schematic version of "Riemannian geometry" we may on a class of affine schemes define the notion "algebraic Levi-Civita connection". On an affine algebraic "Riemannian scheme" $(X,g)$ there is a unique algebraic connection $\nabla$ on the tangent bundle $T_X$ of $X$, that is torsion free and compatible with the "Riemannian metric" $g$. Kunz constructs the Levi-Civita connection $\nabla$ indirectly using the metric $g$, hence his methods differ from the methods in this paper. The methods in this paper constructs all connections on $T_X$ including the connection in Kunz papers (see \cite{kunz0} an \cite{kunz}). I will give explicit examples of connections on vector bundles on affine schemes that are not "Riemannian" in the sense of \cite{kunz0} hence the methods of Kunz do not apply. The papers of Kunz translate some notions and constructions from Riemannian geometry (see \cite{docarmo}) to the language of affine schemes.

Using Theorem \ref{curvatureintro}  I'm able to give explicit formulas for connections on tangent and cotangent bundles on any regular affine hypersurface. The tangent and cotangent bundle
on a regular affine hypersurface are stably trivial hence their Chern classes are zero. The new characteristic class introduced will be non trivial in these cases - I will prove this in the case of the cotangent bundle on the real two sphere. In Corollary \ref{euler} I prove the fundamental class is stronger than the Euler class.
Using explicit calculation I will give explicit formulas for a connection $\nabla$ on the cotangent bundle on the real two sphere where the curvature $R_{\nabla}$ is not a connection of curvature type $f$ where $f$ is a 2-cocycle for the Lie-Rinehart complex. Hence the examples I get give new examples of non flat connections that have not appeared elsewhere in the litterature. I also calculate "by hand" a class of examples that cannot be calculated using computer algebra and Groebner bases.

The new characteristic class $c(E)$ - the \emph{fundamental class of $E$} -  is a class

\[  c(E) \in \Ext^1(L, \End_A(E)) \]

where $\Ext^1(L, \End_A(E))$ is a pointed cohomology set with the structure of a torsor on the abelian group $\operatorname{H}^2(L, Z(\End_A(E))$.   The torsor $\Ext^1(L, \End_A(E))$ classifies equivalence classes of non abelian extensions 

\[  0 \rightarrow \End_A(E) \rightarrow E \rightarrow L \rightarrow 0 \]

of the $A/k$-Lie-Rinehart algebra $L$ by the $A$-Lie algebra $\End_A(E)$. If $E$ is an invertible $A$-module, it follows $\End_A(E) \cong A$, and in this case there is an isomorphism $\Ext^1(L, \End_A(E)) \cong \H^2(L,A)$
where the latter is the 2'nd Lie-Rinehart cohomology of $L$ with values in $A$.  Any  $L$-connection $(E, \nabla)$ gives rise to a non-abelian extension

\[ s(E):  0 \rightarrow \End_A(E) \rightarrow \D(L,E, \nabla) \rightarrow L \rightarrow 0 ,\]

where $\D(L, E, \nabla):=\End_A(E) \oplus L$ and with the following Lie product: For any pair $(f,x),(g,y) \in \D(L,E, \nabla)$, let $\nabla_1$ be the adjoint connection and define the following bracket:

\[  [(f,x), (g,y)]:= ([f,g]+ \nabla_1(x)(g) - \nabla_1(y)(f) + R_{\nabla}(x,y), [x,y] ).\]

This gives $\D(L,E,\nabla)$ the structure of an $A/k$-Lie-Rinehart algebra with a canonical projection map to $L$, making the sequence an "exact sequence" of $A/k$-Lie-Rinehart algebras.
There is a left action 

\[ \sigma: \operatorname{H}^2(L, Z(\End_A(E)) \times \Ext^1(L, \End_A(E)) \rightarrow \Ext^1(L, \End_A(E)) \]

of the underlying abelian group $\operatorname{H}^2(L , Z(\End_A(E))$ on the set $\Ext^1(L, \End_A(E))$. The group $\operatorname{H}^2(L , Z(\End_A(E))$ is defined in terms of the Lie-Rinehart cohomology of the centre $Z(\End_A(E))$ with respect to the adjoint connection on $E$.   In Theorem \ref{mainclass} I prove the following:

\begin{theorem} \label{mainclassintro} Let $(L, \alpha)$ be a $A/k$-Lie-Rinehart algebra where $L$ is projective as $A$-module, and let $(E, \nabla)$ be an arbitrary $L$-connection. The fundamental exact sequence $s(E)$ defines a cohomology class $c(E) \in \Ext^1(L, \End_A(E))$ that is independent of choice of connection $\nabla$. The class $c(E)$ is trivial if and only if $E$ has a flat connection. Let $Z(\End_A(E))$ be the center of the ring of endomorphisms $\End_A(E)$. The adjoint connection $\nabla_1$ induce a flat connection $\nabla'$ on $Z(\End_A(E))$.
It follows there is a group sction

\[ \sigma: \operatorname{H}^2(L, Z(\End_A(E)) \times \Ext^1(L, \End_A(E)) \rightarrow \Ext^1(L, \End_A(E)) \]

making $\Ext^1(L, \End_A(E))$ a torsor on the abelian group $\operatorname{H}^2(L, Z(\End_A(E))$. The group   $\operatorname{H}^2(L, Z(\End_A(E))$  is the the Lie-Rinehart cohomology of $(Z(\End_A(E)), \nabla')$.
\end{theorem}

If $E$ is a finite rank projective $A$-module, there is an isomorphism $Z(\End_A(E)) \cong A$ and it follows $\Ext^1(L, \End_A(E))$ is a torsor on $\H^2(L,A)$.

Note: In the case when $k$ is a field of characteristic $p>0$ and if $k \subseteq L$ is a purely inseparable extension of exponent one, there is an isomorphism of groups 
$\H^2_{res}(\Der_k(L), L) \cong \operatorname{Br}^L_k$, where $\H^2_{res}(\Der_k(L),L)$ is the "restricted Lie-Rinehart cohomology" of the "restricted Lie-Rinehart algebra" $\Der_k(L)$ with values in $L$ and $\operatorname{Br}^L_k$ is the Brauer group. Hence in characteristic $p>0$ it may be there is an action of the Brauer group on the cohomology set $\Ext^1(L, \End_A(E))$. The Brauer group of a field $K$ can also be defined using etale cohomology: There is an isomorphism $\operatorname{Br}(K) \cong \H^2_{et}(\Spec(K)_{et}, \mathbb{G}_m)$ where the latter is the etale cohomology group of the multiplicative group scheme $\mathbb{G}_m$. The etale cohomolgy group $\H^2_{et}(X_{et}, \mathbb{G}_m)$ also classifies "gerbes" on $X$.

If $c(E)$ is trivial it follows all Chern classes $c_i(E)=0$ are zero. Hence the characteristic class $c(E)$ contains information on all higher Chern classes $c_i(E)$ of $E$ and all other characteristic classes of $E$ defined in terms of 
Chern classes or the curvature of a connection. 


In the case when $E:=T_X$ is the tangent bundle we may view the class $c(T_X)$ as an obstruction for $X$ to be "algebraically parallelizable": If the tangent bundle $T_X$ is trivial we must have $c(T_X)$ to be trivial since any trivial bundle has a flat algebraic connection.  The class $c(E)$ is defined for any commutative ring $A$ and any $A$-module $E$ with a connection.

Note: An important and classical - but not so well known - result of Suslin says the following: If the base field $k$ is the field of complex numbers and if $A$ is a finitely generated regular $k$-algebra with $X:=\Spec(A)$ and $E$ a finite rank projective $A$-module that is stably  trivial, then $E$ is trivial. Hence the Theorem implies that the cotangent bundle $\Omega^1_X$ and tangent bundle $T_X$ are trivial. Hence any such $X$ is "algebraically parallelizable" when the base field $k$ is the field of complex numbers. The class introduced in this paper is defined for any commutative ring $A$ over any commutative base ring $k$.






\section{An explicit formula for the curvature of an algebraic connection}

In the following let $A$ be a commutative unital $k$-algebra where $k$ is a commutative unital ring and let $E$ be a left $A$-module. Let $A^n:=A\{e_1,..,e_n\}$ be the free rank $n$ $A$-module on the basis $B:=\{e_1,..,e_n\}$
and let $p: A^n \rightarrow E$ and $s: E \rightarrow A^n$ be $A$-linear maps with $p \circ s = Id_E$. Hence the $A$-module $E$ is a finite rank projective $A$-module. Let $\phi:= s \circ p$. It follows $\phi$ is an idempotent element with $\phi^2 = \phi$ and $Im(\phi) \cong E$. Let $F:=ker(\phi)$. There is an isomorphism of $A$-modules $E\oplus F \cong A^n$. For any map $x\in \End_k(A)$ and any map $\psi \in \End_A(A^n)$ with $[\psi]_B:=[b_{ij}]$ define the new map $x(\psi) \in \End_A(A^n)$ as follows: the endomorphism $x$ acts on the coefficients $b_{ij}$ defining $\psi$ in the basis $B$. We let the endomorphism $x(\psi)$ be defined by the following matrix in the basis $B$:

\[  [x(\psi)]_B :=[x(b_{ij})] .\]

If $x\in \End_k(A)$ is an endomorphism we define the following map: $\rho_B(x): A^n \rightarrow A^n$ is defined as follows: If $u \in A^n$ is an element 

\[ [ u ]_B=
\begin{pmatrix} u_1  \\
 u_2 \\
\vdots \\
  u_n 
\end{pmatrix}
\]

we define 

\[ [ \rho_B(x)(u) ]_B=
\begin{pmatrix} x(u_1)  \\
x( u_2) \\
\vdots \\
  x(u_n) 
\end{pmatrix}.
\]

Hence the $k$-linear endomorphism $\rho_B(x)$ acts via $x$ in each coordinate (when expressed in the basis $B$). We get when $x\in \Der_k(A)$  the following formula for any $a\in A, u \in A^n$:

\[  \rho_B(x)(au)=a\rho_B(x)(u)+ x(a)u .\]

Hence the endomorphism $\rho_B(x) \in \Diff^1(A^n)$ is a first order differential operator.

\begin{proposition}  \label{derivation} Let $\psi_1, \psi_2 \in \End_A(A^n)$ be two endomorphisms with $[\psi_1]_B:=[a_{ij}], [\psi_2]_B:=[b_{ij}]$. It follows for any derivation $x \in \Der_k(A)$ we get equalities of endomorphisms

\[  \rho_B(x)  \circ \psi_1 = \psi_1 \circ \rho_B(x) + x(\psi_1),  x(\psi_1 \circ \psi_2)=x(\psi_1) \circ \psi_1 + \psi_1 \circ x(\psi_2).\]

\end{proposition}
\begin{proof} The proof is an explicit calculation using the basis $B$ and the definitions given above.
\end{proof}

\begin{definition} Let $(L, \alpha)$ be an $A/k$-Lie-Rinehart algebra and let $E$ be a left $A$-module. Let $\psi \in \End_A(E)$. An $(L,\psi)$-connection is an $A$-linear map
$\nabla: L \rightarrow \End_k(E)$ with the property that for any $x \in L, a\in A, e \in E$ it follows

\[ \nabla(x)(ae)=a \nabla(x)(e) + \alpha(x)(a) \psi(e). \]

An \emph{ordinary connection} is an $(L, Id_E)$-connection.

The curvature of an $(L,\psi)$-connection $\nabla$ is the following map:

\[ R_{\nabla}: \wedge^2 L \rightarrow \End_k(E) \]

defined by 

\[ R_{\nabla}(x \wedge y ):=[\nabla(x), \nabla(y)]-\nabla([x,y]).\]

Let $\Conn^{tw}(E)$ be the set of $(L, \psi)$-connections on $E$ for varying endomorphisms $\psi \in \End_A(E)$. An object $(E, \nabla_E, \phi_E) \in \Conn^{tw}(E)$ is a \emph{twisted connection}.
For a fixed endomorphism $\psi \in \End_A(E)$ let $\Conn^{\psi}(E)$ be the set of $(L, \psi)$-connections
$\nabla$. The subset $\Conn^{Id_E}(E)$ is the set  of $(L, Id_E)$-connections $\nabla$.  Given two elements $(E_i, \nabla_i, \phi_i) \in \Conn^{tw}(E_i)$ for $i=1,2$, a morphism of twisted connections $f$ is a map of $A$-modules $f: E_1 \rightarrow E_2$ commuting with the induced action of the connections $\nabla_i$. We get a category - the \emph{category of twisted connections},  denoted $\Conn^{tw}(A)$.
\end{definition}

\begin{example} If $\nabla: L \rightarrow \End_k(E)$ is an ordinary $L$-connection and $\psi \in \End_A(E)$ it follows $\onabla(x):= \nabla(x) \circ \psi$ is an $(L, \psi)$-connection.
If $\psi=Id_E$ and $\nabla$ is an ordinary $L$-connection, it follows the curvature gives a map of $A$-modules

\[ R_{\nabla}: \wedge^2 L \rightarrow \End_A(E).\]

\end{example}

\begin{lemma}  The set $\Conn^{tw}(E)$ is a left $A$-module. Let $G:=\Hom_A(L, \End_A(E))$ be the abelian group of potentials on $E$. There is a right action 

\[ \sigma: \Conn^{tw}(E) \times G \rightarrow \Conn^{tw}(E) \]

defined by $\sigma(\nabla, P):=\nabla + P$. Hence the $A$-module $\Conn^{tw}(E)$ is a torsor on the abelian group $G$. For any fixed element $\psi \in \End_A(E)$ it follows the action $\sigma$ induce a right action of $G$ on $\Conn^{\psi}(E)$. Hence the set of ordinary connections $\Conn^{Id_E}(E)$ is a sub torsor of $\Conn^{tw}(E)$ under $\sigma$.
\end{lemma}
\begin{proof} The proof is straight forward.
\end{proof}

Note: The set of ordinary connections $\Conn^{Id_E}(E)$ is not an $A$-module.

Let $E$ be a left $A$-module and let $p:A^n \rightarrow E, s:E \rightarrow A^n$ be $A$-linear maps with $p \circ s =Id_E, \phi:= s \circ p$. Let $x \in L:=\Der_k(A)$ and define the following map

\[ \nabla_E: L \rightarrow \End_k(E) \]

by 

\[ \nabla_E(x)(e):= p \circ \rho_B(x) \circ s(e). \]

It follows $\nabla_E$ is an $A$-linear map with the property that for any $a \in A, e \in E$ we get 

\[ \nabla_E(x)(ae)=a\nabla(x)(e) + x(a)e.\]

Hence the map $\nabla_E$ is an $L$-connection. Define the map $\onabla_E: L \rightarrow \End_k(A^n)$ by

\[  \onabla_E(x)(u):= s \circ \nabla(x) \circ p(u).\]

For any $a\in A, u \in A^n$ we get

\[  \onabla_E(x)(au)=a\onabla_E(x)(u)+ x(a)\phi(u).\]

Hence the map $\onabla_E$ is an  $(L,\phi)$-connection on $A^n$.  For any $x \in L$ there is an equality

\[ \onabla_E(x)= \phi \circ \rho_B(x) \circ \phi. \]

There are canonical $A$-linear maps $\pi_k: \End_k(A^n) \rightarrow \End_k(E)$ defined by $\pi_k( z):= p \circ z \circ s$  and $\gamma_k: \End_k(E) \rightarrow \End_k(A^n)$ 
defined by $\gamma_k(z):= s \circ z \circ p$, and there is a commutative diagram

\[
\diagram          &     \End_k(A^n) \dto^{\pi_k}    \\
                      L \urto^{\onabla_E} \rto^{\nabla_E}  &   \End_k(E)    
\enddiagram.
\]

There are maps of $A$-modules $\pi_A: \End_A(A^n) \rightarrow \End_A(E)$ defined by $\pi_A( z):= p \circ z \circ s$  and $\gamma_A: \End_A(E) \rightarrow \End_A(A^n)$ 
defined by $\gamma_A(z):= s \circ z \circ p$.

\begin{proposition} Let $\nabla_E$ be the connection defined above and let $\onabla_E:= \gamma_k \circ \nabla_E$. It follows $\onabla_E$ is an $(L, \phi)$-connection on $A^n$. The curvature $R_{\onabla_E}$ gives an $A$-linear map

\[ R_{\onabla_E}: \wedge^2 L \rightarrow \End_A(A^n) \]

and an equality $R_{\nabla_E}=\pi_A \circ R_{\onabla_e}$. Hence we recover the curvature of $\nabla_E$ from the curvature of $\onabla_E$. Let $\nabla_1: L \rightarrow \End_k(E)$ be any connection. It follows there is a "potential" $P \in \Hom_A(L, \End_A(E))$ with $\nabla_1 = \nabla_E + P$ and an equaltiy

\[  R_{\nabla_1}(x,y)= R_{\nabla_E}(x,y)+ d^1_{\onabla_E}(P)(x,y) + [P(x), P(y)], \]

hence the curvature of $\nabla_1$ may be recovered from the curvature of $\nabla_E$ and the potential $P$.
\end{proposition}
\begin{proof} A calculation shows that the curvature $R_{\onabla_E}(x,y) \in \End_A(A^n)$ and  that $\pi_A \circ R_{\onabla_E}=R_{\nabla_E}$. The rest of the Proposition follows directly from the definitions.
\end{proof}

There is a commutative diagram

\[
\diagram          &     \End_A(A^n) \dto^{\pi_A}    \\
                      \wedge^2 L \urto^{R_{\onabla_E}} \rto^{R_{\nabla_E}}  &   \End_A(E)    
\enddiagram.
\]

Hence the $(L,\phi)$-connection $\onabla_E$ is a "lifting" of the $L$-connection $\nabla_E$ to $A^n$ and we recover the curvature $R_{\nabla_E}$ from $R_{\onabla_E}$. The aim is to calculate the curvature of $\nabla_E$ in terms of the curvature of the $(L,\phi)$-connection $\onabla_E$.
The lifted $(L,\phi)$-connection $\onabla_E$ is a connection on a free module $A^n$, and since $\phi$ is an idempotent element it follows the curvature $R_{\onabla_E}(x,y)$ is an $A$-linear map for any two elements $x,y \in L$. It follows we may express the curvature of $\onabla_E(x,y)$ in a given basis. The aim is to give an explicit formula for the curvature of $\onabla_E(x,y)$ in a basis.

\begin{example} An exact sequence of twisted connections. \end{example}

Let 

\[  0 \rightarrow F \rightarrow A^n \rightarrow^p E \rightarrow 0 \]

be the exact sequence induced by the surjective map $p$.  Let $\nabla_0: L \rightarrow \End_k(F)$ be defined by $\nabla_0(\partial)(u)=0$ for all $\partial \in L, u \in F$ and let $\phi_0:=0$ be the zero endomorphism
of $F$. It follows $\nabla_0$ is an $(L, \phi_0)$-connection on $F$. There is an exact sequence of twisted connections  in $\Conn^{tw}(A)$

\[
\diagram    0 \rto      &   F \rto \dto^{\nabla_0(\partial)}  &        A^n \rto \dto^{\onabla_E(\partial)}   & E \rto \dto^{\nabla_E(\partial)}  & 0    \\
                  0 \rto    &   F \rto                                            &        A^n \rto                                & E \rto                               & 0
\enddiagram.
\]

 The above construction is functorial in the following sense: Let $(L,\alpha)$ be an $A/k$-Lie-Rinehart algebra and let $\Conn^s(L)$ denote the following category. The objects are tuples
$(E, \nabla_E, A^{n(E)}, p_E, s_E)$ where $(E, \nabla_E)$ is an $L$-connection, $n(E) \geq $ an integer and $p_E: A^{n(E)} \rightarrow E, s_E: E \rightarrow A^{n(E)}$ maps of $A$-modules with $p_E \circ s_E = Id_E$
a split surjection.  A map 

\[  f: (E, \nabla_E, A^{n(E)}, p_E, s_E) \rightarrow (F, \nabla_F, A^{n(F)}, p_F, s_F) \]

is a map of $A$-modules $f: E \rightarrow F$ that is a map of connections. 


\begin{proposition}\label{functor}  There is a functor

\[ G: \Conn^s(L) \rightarrow \Conn^{tw}(A) \]

defined as follows: Define for any object $(E, \nabla_E, A^{n(E)}, p_E, s_E)$ the map $\onabla_E: L \rightarrow \End_k(A^{n(E)})$ by $\onabla(x)(u):= s_E \circ \nabla_E(x) \circ p_E$.  Define $G$ as follows:

\[  G(E, \nabla_E, A^{n(E)}, p_E, s_E):=(A^{n(E)}, \onabla_E).\]

For any map of connections
$f: E\rightarrow F$ define $G(f):= \tilde{f}:=s_F \circ f \circ p_E$.  It follows we get a map of connections $\tilde{f}:(A^{n(E)}, \onabla_E) \rightarrow (A^{n(F)}, \onabla_F)$.
\end{proposition}
\begin{proof} The proof follows from the above definitions.
\end{proof}

Let $x,y \in L$ be two derivations and consider the curvature $R_{\onabla_E}(x,y)$. Using Proposition \ref{derivation} we get the following:

\[ R_{\onabla_E}(x,y):= \onabla_E(x)\onabla_E(y)-\onabla_E(y)\onabla_E(x)-\onabla_E([x,y] ) :=\]

\[ \phi \circ \rho_B(x) \circ \phi \circ   \rho_B(y) \circ \phi -     \phi \circ \rho_B(y) \circ \phi  \circ \rho_B(x) \circ \phi  - \phi \circ \rho([x,y]) \circ \phi  = \]

\[   \phi \circ \rho_B(x) \circ \phi ( \phi \circ \rho_B(y)  + y(\phi)) \]

\[ -  \phi \circ \rho_B(y) \circ \phi ( \phi \circ \rho_B(x)  + x(\phi)) \]

\[  - \phi \circ \rho_B([x,y]) - \phi \circ [x,y](\phi) = \]

\[   \phi \circ \rho_B(x) \circ  \phi \circ \rho_B(y)   + \phi \circ \rho_B(x) \circ \phi \circ y(\phi)) \]

\[ -  \phi \circ \rho_B(y) \circ  \phi \circ \rho_B(x)  -  \phi \circ \rho_B(y) \circ \phi \circ     x(\phi)) \]

\[  - \phi \circ \rho_B([x,y]) - \phi \circ [x,y](\phi) = \]

\[   \phi \circ ( \phi \circ \rho_B(x) + x(\phi)) \circ \rho_B(y) + \phi \circ ( \phi \circ y(\phi) \circ \rho_B(x) + x( \phi \circ y(\phi))) \]

\[-   \phi \circ ( \phi \circ \rho_B(y) + y(\phi)) \circ \rho_B(x) - \phi \circ ( \phi \circ x(\phi) \circ \rho_B(y) + y( \phi \circ x(\phi))) \]

\[ - \phi \circ \rho_B([x,y]) - \phi \circ [x,y](\phi) =  \]

\[ \phi (\phi \circ \rho_B(x) + x(\phi) ) \circ \rho_B(y) + \phi \circ (\phi \circ y(\phi) \circ \rho_B(x) + x(\phi \circ y(\phi)))  \]

\[-  \phi (\phi \circ \rho_B(y) + y(\phi) ) \circ \rho_B(x) - \phi \circ (\phi \circ x(\phi) \circ \rho_B(y) + y(\phi \circ x(\phi)))  \]

\[ - \phi \circ \rho_B([x,y]) - \phi \circ [x,y](\phi) =  \]

\[  \phi \circ \rho_B(x) \circ \rho_B(y) + \phi \circ x(\phi) \circ \rho_B(y) + \phi \circ y(\phi) \circ \rho_B(x) + \phi \circ (x(\phi) \circ y(\phi) + \phi \circ xy(\phi) ) \]

\[ -  \phi \circ \rho_B(y) \circ \rho_B(x) + \phi \circ y(\phi) \circ \rho_B(x) + \phi \circ x(\phi) \circ \rho_B(y) + \phi \circ (y(\phi) \circ x(\phi) + \phi \circ yx(\phi) ) \]

\[ - \phi \circ \rho_B([x,y]) - \phi \circ [x,y](\phi) =  \]

\[ \phi \circ [\rho_B(x), \rho_B(y)] - \phi \circ \rho_B([x,y]) + \phi \circ (x(\phi ) \circ y(\phi) - y(\phi) \circ x(\phi) )  \]

\[ + \phi \circ (  xy(\phi) - yx(\phi) ) - \phi \circ [x,y](\phi) =\]

\[  \phi \circ [x(\phi), y(\phi) ].\]

We may prove the main Theorem of this section.

\begin{theorem}\label{curvature}  Let $E$ be a finite rank projective $A$-module and let $p: A^n \rightarrow E, s:E \rightarrow A^n$ be $A$-linear maps with $p \circ s =Id_E$. Let $\phi:=s \circ p$ be the idempotent element defining $E$. Let $B:=\{e_1,..,e_n\}$ be a basis for $A^n$ and let $\nabla_E: \Der_k(A) \rightarrow \End_k(E)$ be the connection defined by $\nabla_E(x)(e):= p \circ \rho_B(x) \circ s(e)$. Let $R_{\nabla_E}$ be the curvature of $\nabla_E$. Let $\onabla_E: \Der_k(A) \rightarrow \End_k(A^n)$ be the induced $(L, \phi)$-connection from Proposition \ref{functor}  with curvature $R_{\onabla_E}$. For any elements $x,y \in \Der_k(A)$ it follows 

\[ R_{\onabla_E}(x,y) = \phi \circ [x(\phi), y(\phi)] .\]

We recover the curvature $R_{\nabla_E}(x,y)$ via the formula

\[  R_{\nabla_E}(x,y)= \pi_A( R_{\onabla_E}(x,y)).\]

\end{theorem}
\begin{proof} Given any two derivations $x,y \in \Der_k(A)$ it follows there is by the calculation above an equality of endomorphisms of $A^n$:

\[  R_{\onabla_E}(x,y)= \phi \circ [x(\phi), y(\phi)] .\]

The Theorem follows.

\end{proof}

If we express the endomorphisms from Theorem \ref{curvature} in the  basis $B$ for the free $A$-module $A^n$, we get an explicit formula for the curvature of the lifted $(L, \psi)$-connection $\onabla_E$ in the basis $B$:

\begin{corollary}  \label{maincorollary} Assume $A^n:=A\{e_1,..,e_n\}$ where $B:=\{e_i\}$ is a basis as $A$-module, and let $[\phi]_B:=[a_{ij}]$ be the matrix of the idempotent element $\phi$ in the basis $B$, where $\phi$ is an idemopotent element defining the $A$-module $E$.
It follows the curvature $R_{\onabla_E}$ from  Theorem \ref{curvature} is as follows: Let  $x,y \in \Der_k(A)$  be derivations and let $[x(a_{ij})]$ be the matrix with $x(a_{ij})$ in the $(i,j)$-spot. Let $[R_{\onabla_E}(x,y)]_B$ be the matrix of the curvature
$R_{\onabla_E}$ in the basis $B$.  It follows there is an equality of matrices

\[  [R_{\onabla_E}(x,y)]_B =  [ a_{ij}]( [x(a_{ij})] [y(a_{ij})] - [y(a_{ij}) ][x(a_{ij})] ).\]

We get an explicit formula for curvature $R_{\nabla_E}$ of the connection $\nabla_E$  using Theorem \ref{curvature} and the formula $R_{\nabla_E}=\pi_A \circ  R_{\onabla_E}$.
\end{corollary}
\begin{proof} The Corollary follows from Theorem \ref{curvature}.
\end{proof}

Let $\nabla_1$ be any connection on $E$. It follows there is a potential $P$ an an equality $\nabla_1=\nabla_E + P$ and the curvature $R_{\nabla_1}$  satisfies the following formula for any $x,y \in L$:

\[  R_{\nabla_1}(x,y)=R_{\nabla_E}(x,y) + d^1_{\onabla_E}(P)(x,y) + [P(x), P(y)] .\]

Hence we get a formula for the curvature of $\nabla_1$  in terms of the curvature of $\nabla_E$ and the potential $P$. Hence Corollary \ref{maincorollary} gives an explicit formula for the curvature
of any connection on $E$ in terms of an idempotent element $\phi$ and a potential $P$.

Since the formula $R_{\onabla_E}(x,y):=\phi \circ [x(\phi), y(\phi)]$ from Theorem \ref{curvature} involves the Lie product $[x(\phi), y(\phi)]$ between the two endomorphisms $x(\phi)$ and $y(\phi)$, it follows "most" connections on $E$
are non flat. It is a non trivial problem to determine if $E$ has a flat algebraic connection.

\section{Explicit formulas for algebraic connections on regular hypersurfaces}

In this section I will give explicit examples illustrating Theorem \ref{curvature} in the case when $X$ is an affine regular complete intersection and $E$ is the tangent or contangent bundle on $X$.
The tangent bundle (and cotangent bundle) on $X$ is stably trivial hence all its Chern classes are trivial. It is a theorem of Suslin that the cotangent and tangent bundle of any regular affine complete intersection over the complex numbers is trivial. Hence a large class of affine varieties have tangent and cotangent bundles with vanishing Chern classes. For this reason it is interesting to search for new invariants that are non trivial for such vector bundles.  In this section I will study a characteristic class $c(E)$ introduced in another paper (see \cite{maa}) and  I will prove that the class $c(T_X)$ and $c(\Omega_X)$ is non trivial in many cases. I will also give explicit examples 
of non flat connections that are not of curvature type $f$ where $f$ is a 2-cocycle for the Lie-Rinehart complex. Hence the connections I construct give new non trivial examples of non flat connections that have not appeared
elsewhere in the litterature.

Note: Kunz constructs in \cite{kunz0} a connection $\nabla$ on the tangent bundle $T_X$ of an affine Riemannian scheme $(X,g)$ over a field $k$ of characteristic different from $2$ using a "metric" $g$ on the tangent bundle. The metric is a symmetric bilinear non-degenerate form on $T_X$

\[  g(-,-): T_X \oplus T_X \rightarrow A .\]

Kunz defines a "Levi-Civita connection " $\nabla$ implicitly using the following formula: Let $x,y,z \in T_X$ be vector fields and define 

\[  g(z, \nabla(y)(x)) = \frac{1}{2}( x(g(y,z)) + y(g(z,x)) - z(g(x,y))    \]

\[ -  g([x,z],y)- g([y,z], x) - g([x,y],z)) .\]

Then he proves there is a unique connection $\nabla$ satisfying this condition for all $x,y,z$, that is torsion free and compatible with the metric $g$. The connection $\nabla$ is the "algebraic Levi-Civita connection" on $T_X$. He also proves that a regular hypersurface $X$ defined by a polynomial 
$f (x_1,..,x_n) \in k[x_1,..,x_n]$ has a Riemannian metric on $T_X$ if and only if the following holds: $  \sum f_{x_i}^2$ is a unit in the coordinate ring $A:=k[x_i]/(f)$. Here $f_{x_i}$ is the partial derivative of $f$ with respect to $x_i$. Hence not all regular
hypersurfaces has a Riemannian metric in the sense of Kunz' preprint \cite{kunz0}. The main aim of this paper is to study the characteristic class $c(E)$ for arbitrary $E$ on an arbitrary affine scheme $X$ and for this reason we must introduce more general techniques.

One wants to relate the class $c(E)$ to classical notions uch as the \emph{Euler characteristic} and to give schematic proofs of classical formulas from differential geometry. For an affine regular surface $S$ with a Riemannian metric $g$ one wants to calculate the Euler characteristic  $\chi(S)$ using the class $c(T_S) \in \Ext^1(T_S, \End(T_S))$. The Euler characteristic $\chi(S)$ cannot be calculated using the Chern classes $c_i(T_S)$ since they are all zero. If such a relation exists this would give an algebraic
definition of the Euler characteristic of $S$. The Euler characteristic is a topological invariant defined in terms of the Betti cohomology of $S$ when $S$ is a real smooth surface.

Note: By \cite{kunz} Theorem 9.6 it follows for any $k$-algebra $A$ which is a finitely generated regular integral domain such that the extension $k \subseteq K(A)$ is separable, where $K(A)$ is the quotient field
the following holds: If $A$ is a global complete intersection and if $dim(A)=1$ it follows $\Omega^1_A$ is a free finite rank $A$-module. In this case it follows the module $\Omega^1_{A/k}$ has a flat algebraic connection and hence $c(E)$ is trivial.  If $C:=\Spec(A) \subseteq \mathbb{A}^2_k$ is a plane affine curve one easily checks via explicit calculations there is a trivialization of $\Omega^1_{A/k}$.

In the following let $A$ be any commutative unital $k$-algebra with $k$ a commutative unital ring, and let $E$ be a left $A$-module with a connection $\nabla_E: L \rightarrow \End_k(E)$ where $(L, \alpha)$ is an $A/k$-Lie-Rinehart algebra.
Let $\onabla_E: L \rightarrow \End_k(\End_A(E))$ be the induced connection defined by $\onabla_E(x)(f):=[\nabla_E(x), f]$ where $[,]$ is the Lie product on $\End_k(E)$. Let $\End(E,L):=\End_A(E) \oplus L$ with the following product:

\[  [(f,x), g,y)]:= ([f,g] + \onabla_E(x)(g)- \onabla_E(y)(f)+ R_{\nabla_E}(x,y), [x,y])  .\]

Let $\alpha_E: \End(E,L) \rightarrow L$ be the canonical projection map. It follows the pair $(\End(E,L), \alpha_E)$ is an $A/k$-Lie-Rinehart algebra. 
There is an $A$-linear  map 

\[ \nabla_1: \End(E,L) \rightarrow \End_k(E) \]

 defined by $\nabla_1(f,x)(g):=[f+ \nabla_E(x), g]$. 

\begin{lemma} The map $\nabla_1:\End(E,L) \rightarrow \End_k(E)$ is a flat $\End(E,L)$-connection. 
\end{lemma}
\begin{proof} The proof is an elementary calculation with connections.
\end{proof}

Hence for an arbitrary $L$-connection $(E,\nabla_E)$ there is always an extention of $A/k$-Lie-Rinehart algebras $p: \End(E,L) \rightarrow L$ and an induced connection $\nabla_1$ on $\End(E,L)$ that is flat. 

\begin{definition} \label{curvaturetype} Let $(E, \nabla_E)$ be an $L$-connection. Let $ \rho \in \End_A(E)$ be an endomorphism and let $f \in Z^2(L,A)$ be a 2-cocycle for the Lie-Rinehar complex of $L$.
We say the connection $R_{\nabla_E}$ has \emph{curvature type $(f, \rho)$} iff for all $x,y \in L$ and $e \in E$ it follows 

\[  R_{\nabla}(x,y)(e)=f(x,y)\rho(e).\]

If $\rho=Id_E$ we say the connection $\nabla_E$ has \emph{curvature type $f$}.
\end{definition}

Note: If the 2-cocycle $f$ is nonzero, it follows a connection of the form from Definition \ref{curvaturetype} is a special case of a non flat connection.

Note: Let $A$ be an integral domain and let $\nabla_E$ be a connection on $E$ with the property that there is an endomorphism $\rho \in \End_A(E)$ where the curvature $R_{\nabla_E}(x,y)$ is a multiple of $\rho$ for any $x,y \in L$. 
It follows there is a 2-cocycle $f \in Z^2(L,A)$ and an equality $R_{\nabla_E}(x,y)(e)=f(x,y)\rho(e)$. Hence $\nabla_E$ has curvature type $(f,\rho)$.

For simplicity let us study the case of an affine regular surface defined by a polynomial $f(x,y,z)$ over an arbitrary commutative unital ring $k$.
Let $x,y,z$ be independent variables over $k$ and let $f(x,y,z) \in k[x,y,z]$ be a polynomial and let $A:=k[x,y,z]/(f)$. Assume the jacobian ideal $J(f):=(f_x, f_y,f_z)=(1)$ is the unit ideal in $A$. It follows there are elements $a,b,c, h \in k[x,y,z]$ and an equality in $k[x,y,z]$

\[  af_x+ bf_y + cf_z = 1 + hf .\]

Let $B:=\{ dx,dy,dz \}$ be a basis for the free $A$-module $F:=A\{dx,dy,dz\}$ of rank $3$ and Let $Df:=f_xd_x f_ydy + f_zdz \in F$. By definition there is an exact equence of $A$-modules

\[  0 \rightarrow ADf \rightarrow F \rightarrow \Omega  \rightarrow 0, \]

where $\Omega:=\Omega^1_{A/k}$ is the module of Kahler differentials of $A/k$.  There is an $A$-linear splitting $s$ of the canonical projection map $p: F \rightarrow \Omega$ defined as follows:

\[ s(\overline{dx}) :=dx - aDf, s(\overline{dy}):=dy - b Df, s(\overline{dz}):=dz-cDf.\]

It follows the endomorphism $\psi := s \circ p \in \End_A(F)$ is an idempotent element defining the module $\Omega$: There is an isomorphism $Im(\psi) \cong \Omega$. In the basis $B$ we get the following matrix defining $\psi$:

\[ 
[\psi]_B = [a_{ij}] :=
\begin{pmatrix}   1-af_x   &  -bf_x   &   -cf_x  \\
                             -af_y   &  1-bf_y    & -cf_y \\
                            -af_z  & -bf_z  & 1-cf_z   
\end{pmatrix}.
\]

By the results in the previous section we get a connection $\nabla_{\Omega}$ and an induced $(L, \psi)$-connection $\onabla_{\Omega}$ and two commutative diagrams

\[
\diagram          &     \End_k(F) \dto^{\pi_k}    \\
                      L \urto^{\onabla_{\Omega}} \rto^{\nabla_{\Omega}}  &   \End_k(\Omega)    
\enddiagram.
\]

There are maps of $A$-modules $\pi_A: \End_A(F) \rightarrow \End_A(\Omega)$ defined by $\pi_A( z):= p \circ z \circ s$  and $\gamma_A: \End_A(\Omega) \rightarrow \End_A(F)$ 
defined by $\gamma_A(z):= s \circ z \circ p$, and there is a commutative diagram

\[
\diagram          &     \End_A(F) \dto^{\pi_A}    \\
                      \wedge^2 L \urto^{R_{\onabla_{\Omega}}} \rto^{R_{\nabla_{\Omega}}}  &   \End_A(\Omega)    
\enddiagram.
\]

If $L:=\Der_k(A)$ and $x,y \in L$ it follows from Corollary \ref{maincorollary} there is an equality of matrices

\[   [R_{\onabla_{\Omega}}(x,y)]_B = [a_{ij}]( [x(a_{ij})][y(a_{ij})] - [y(a_{ij})][x(a_{ij})] ) , \]

where we let the endomorphisms $x$ and $y$ act on the coefficients $a_{ij}$ of the matrix $[\psi]_B$.  By the previous section there is an equality

\[ R_{\nabla_{\Omega}} = \pi_A \circ R_{\onabla_{\Omega}}, \]

since the curvature $R_{\onabla_{\Omega}}$  has been factored through the free module $F$ via the curvature $R_{\onabla_{\Omega}}$ of the $(L, \psi)$-connection $\onabla_{\Omega}$.
 In the following I will calculate the curvature in the case of the 2-sphere.

\begin{example} A non flat connection on the cotangent bundle on the real two sphere.\end{example}

Let $f:=x^2+y^2+z^2-1 \in k[x,y,z]$ and let $A:=k[x,y,z]/(f)$. Let $Df:=xdx+ydy+zdz\in A\{dx,dy,dz\}$ and let $\Omega:=A\{dx,dy,dz\}/ADf$ be the module of differentials of $A/k$.  Let

\[ 
\rho := 
\begin{pmatrix}   0    &   z   &   -y  \\
                             -z   &   0     &  x \\
                            y   & -x  &  0    
\end{pmatrix}.
\]

And let $M:= \rho \circ \rho^{tr}:=-\rho^2$. It follows $M$ defines an idempotent element $M \in \End_A(A^3)$ with $Im(M) \cong \Omega$.  The endomorphism $\rho \in \operatorname{GL}_A(\Omega)$ and it follows 
$\rho \circ \rho^{tr} = \rho^{tr} \circ \rho = Id_{\Omega}$. We may use the results in the previous section to calculate the curvature $R_{\nabla}$
of the connection $\nabla$ corresponding to the idempotent endomorphism $M$. Let 

\[  D_1:=y \partial_x - x\partial_y, D_2:=z\partial_x - x\partial_z, D_3:=z\partial_y - y \partial_z.\]

It follow the derivations $D_1,D_2, D_3$ generate the module $\Der_k(A)$ as a left $A$-module and we want to calculate the curvature $R_{\nabla}(D_i, D_j)$ using the formula from the previous section.
We get the following theorem:

\begin{theorem}  \label{maincalc} Let $\nabla$ be the connection on $\Omega$ corresponding to the idempotent element $M$. There are equalities   

\[ R_{\nabla}(D_1,D_2)= x\rho,  R_{\nabla}(D_1,D_3)= y\rho,  R_{\nabla}(D_2,D_3)= z\rho \]

in $\End_A(\Omega)$. Hence the connection $\nabla$ has curvature type $(f, \rho)$ where $f$ is a 2-cocycle and $\rho$ is the automorphism of $\Omega$ defined above. 
\end{theorem}
\begin{proof} The proof follows from an explicit calculation using Corollary \ref{maincorollary}.
\end{proof}

The automorphism $\rho$ from Theorem \ref{maincalc} is a non trivial skew symmetric automorphism of $\Omega$ and gives the first example of a naturally ocurring non flat connection 
that is not of curvature type $f$ for a 2-cocycle $f$. The connection $\nabla$ is not an algebraic Levi-Civita connection in the sense of the construction of Kunz.

We may calculate the curvature $R_{\nabla}$ as a map:

\begin{corollary} There is an isomorphism $\wedge^2 \Der_k(A) \cong A w$ where $w:= zd_1 \wedge d_2 - y d_1 \wedge d_3 + xd_2\wedge d_3$  and where $d_1:=D_3, d_2:=-D_2, d_3:=D_3$. Hence the curvature map

\[ R_{\nabla} \in \Hom_A(\wedge^2 \Der_k(A), \End_A(\Omega)) \]

is the following map: $R_{\nabla}(aw):= a(2xz-y^2) \rho$ for $a\in A$.
\end{corollary}
\begin{proof} The proof follows from an explicit calculation.
\end{proof}

\begin{corollary} Let $k$ be the real number field and let $S.=\Spec(A)$ be the real two sphere. The class $c(\Omega)$ is non trivial.
\end{corollary}
\begin{proof} If $\Omega$ had a flat algebraic connection this would imply that $\Omega$ is trivial, which is a contradiction. Hence $c(\Omega)$ is non trivial.
\end{proof}

Note: Since the curvature $R_{\nabla}$ of the connection from Theorem \ref{maincalc} is given by a skew symmetric automorphism $\rho \in \operatorname{GL}_A(\Omega)$ we may ask for a construction of the "Pfaffian"
of the curvature $R_{\nabla}$ similar to the construction used in differential geometry and Chern-Weil theory. Using the Pfaffian we may ask for a construction of the Euler class and the Euler characteristic of the real 2-sphere $S^2_{\mathbb{R}}$ using  the connection from Theorem \ref{maincalc}. This would give a purely algebraic construction of the topological Euler characteristic of the real two sphere.

\begin{corollary} \label{euler}  Let $k$ be the real number field and let $S^n$ be an odd sphere over $k$ with $n \neq 7$. It follows $c(TS^n)$ is non trivial. The Euler class $e(TS^n)$ is trivial hence the fundamental class is stronger than the Euler class.
\end{corollary}
\begin{proof} It follows from  \cite{kobayashi},  Proposition 9.9.2 that $TS^n$ does not have a flat connection hence $c(TS^n)$ is non trivial. It is well known the tangent bundle $TS^n$ has a nowhere vanishing vector field
hence $e(TS^n)$ is trivial.
\end{proof}

\begin{example} A class of interesting  differential forms on the two-sphere. \end{example}

In this Example I construct a countably infinite class of interesting differential 1-forms $\omega_n \in \Omega:=\Omega^1_{A/k}$ where $k$ is a field of characteristic zero and $A:=k[x,y,z]/(f)$
and $f:=x^2+y^2+z^1-1$. There is the algebraic deRham complex

\[ 0 \rightarrow A \rightarrow^{d^0} \Omega^1 \rightarrow^{d^1} \wedge^2 \Omega \rightarrow 0\]

and by definition $\H^1_{dR}(A):= ker(d^1)/Im(d^0)$. 


Let for $n \geq 1$

\[ \omega_n:=(x^2+y^2)^n(x(zdx)+ y(zdy)) \]

where by definition

\[ \Omega:=A\{dx,dy,dx\}/Adf \]

where $df:=xdx+ydy+ zdz$.  This formula follows from Matsumuras book \cite{matsumura}. Let $w:=zdx \wedge dy -ydx \wedge dz + x dy \wedge dz$. It follows $\wedge^2 \Omega \cong Aw$ is the free $A$-module
of rank one on the element $w$. Let $\omega:=x(zdx) + y(zdy)$ and

\[ \omega_n:=(x^2+y^2)^n \omega.\]

It follows

\[ d^1(\omega):=d(xz)\wedge dx + d(yz) \wedge dy = xdz \wedge dx + y dz \wedge dy = (xy)\omega -(yx) \omega =0\]

in $\wedge^2 \Omega$, hence $d^1(\omega)=0$. We get

\[  d^1(\omega_n):= d((x^2+y^2)^n) \wedge \omega = \]

\[ 2n(x^2+y^2)^{n-1}(xyz(dx \wedge dy) + xyz(dy \wedge dx) =0 \]

hence for  all $n \geq 1$ it follows $d^1(\omega_n)=0$ and we get for every $n \geq 0$ a cohomology class 

\[  \overline{\omega_n} \in \H^1_{DR}(S^2).\]

In the next section I will study this construction.


\begin{example}  The Cancellation problem for stably trivial vector bundles. \end{example}

Let $A$ be a commutative unital ring with $S:=\Spec(A)$ and let $Y:=Spec(B)$ be an affine scheme over $S$. The \emph{generalized cancellation problem} asks the following:If there is an isomorphism (of schemes over $S$)

\[  \mathbb{A}^m_S \times_S Y \cong \mathbb{A}^{m+l}_S \]

does this imply there is an isomorphism $Y \cong \mathbb{A}^l_S$ over $S$?  The \emph{classical cancellation problem} asks this question for $A:=k$ a field (see \cite{kraft}). As an example, let $E$ ba a stably trivial $A$-module. This means there is a split exact sequence of $A$-modules

\[ 0 \rightarrow A^m \rightarrow A^n \rightarrow E \rightarrow 0 \]

implying an isomorphism (of schemes over $S$)

\[ \mathbb{A}^m_S \times_S \mathbb{V}(E^*) \cong \mathbb{A}^n_S .\]

For a regular hypersurface it follows the cotangent bundle and tangent bundle are stably trivial.  If $R$ is a fixed commutative base ring and if $f(x_1,..,x_n) \in A:=R[x_1,..,x_n]$ is a polynomial where the jacobian ideal $J(f)=(1)$ is the unit ideal, we get a split exact sequence

\[ 0 \rightarrow Adf \rightarrow A^n  \rightarrow  \Omega^1_{A/R} \rightarrow 0 \]

and an isomorphism

\[  \mathbb{A}^1_S\times_S \mathbb{V}((\Omega^1_{A/R})^*) \cong \mathbb{A}^n_S .\]

A similar result holds for the dual of $\Omega^1_{A/R}$ and the tangent bundle: There is an isomorphism of schemes over $S$:

\[  \mathbb{A}^1_S\times_S \mathbb{V}((T_{A/R})^*) \cong \mathbb{A}^n_S .\]


The invariant introduced in this paper is by this section "computable" and may have applications to the study of this classical problem.
For any stably trivial vector bundle $E$ we may ask this question and if the class $c(E)$  is non trivial we get a counter example to the claim.

\begin{example} A class of regular $3$-folds studied by Russel and the Makar-Limanov invariant. \end{example}

Let $f:=x(1+xy)+z^3+t^2 \in k[x,y,z,t]$ with $k$ the complex number field and let $Y:=Z(f) \subseteq \mathbb{A}^4_k$ be the hypersurface defined by $f$ with coordinate ring $k[Y]:=k[x,y,z,t]/(f)$. As mentioned in \cite{kraft} the following holds:
$Y$ is a smooth hypersurface and diffeomorphic (as reals smooth manifold) to $k^3$. The algebraic hypersurface $Y$ is not isomorphic to $\mathbb{A}^3_k$ as affine scheme over $k$. Since $Y$ is a regular hypersurface it follows its tangent bundle $T_Y$
is stably trivial: There is an isomorphism $T_Y \oplus k[Y] \cong k[Y]^4$. Hence $c_i(T_Y)=0$ have all Chern classes equal to zero. The fundamental class $c(T_Y)$ may be non zero and it may be it can be useful in distinguishing $Y$
from affine $3$-space.   Affine $3$-space have trivial fundamental class and the fundamental class is an isomorphism invariant.

The \emph{Makar-Limanov invariant} of an affine variety $Y$ with coordinate ring $k[Y]$, is an invariant $ML(Y) \subseteq k[Y]$, which is a sub-$k$-algebra of $k[Y]$ defined in terms of locally nilpotent vector fields on $Y$.
The invariaint $ML(Y)$ can be used to distinguish $Y$ from affine $3$-space. The fundamental class $c(T_Y)$ is similar to classical characteristic classes such as the Chern class and Euler class and it may be the fundamental class
can be used to give a systematic approach to the study of the cancellation problem. It is an open problem (according to \cite{kraft}) to decide if all vector bundles on the affine varieties introduced by Russel are trivial. Here the
fundamental class can be useful since it measures if a vector bundle is trivial.

\section{The algebraic deRham cohomology of the two-sphere}

In this section I use the methods developed in the previous section to calculate the algebraic deRham cohomology groups of the two-sphere $S^2$ over any field $k$ of characteristic zero. 
I prove that the cohomolggy group $\H^1_{DR}(S^2)$ is countably infinite dimensional contradicting a classical result from the 1960s relating algebraic deRham cohomology to singular cohomology. 

Let $k$ be any field of characteristic zero and let $f(x,y,z):=z^2-p(x,y)$ with $p(x,y):=1-x^2-y^2$. Let $A:=k[x,y,z]/(f)$ and $S^2:=\Spec(A)$. There is an isomorphism $A \cong k[x,y]\{1, z\}$ and $A$ is a free rank two 
$k[x,y]$-module on the elements $1,z$. Let $\Omega:=\Omega^1_{A/k}$ be the module of Kahler differentials of $A/k$ and let 

\[ 0 \rightarrow A \rightarrow^{d^0} \Omega \rightarrow^{d^1} \wedge^2 \Omega \rightarrow 0 \]

be the algebraic deRham complex of $A/k$. The aim of this section is to calculate the cohomology groups $\H^i_{DR}(S^2):=ker(d^i)/Im(d^{i-1})$ of this complex.

\begin{lemma} There is an isomorphism $\Omega \cong k[x,y]\{dx, zdx, dy, zdy,dz\}$ hence $\Omega$ is a free $k[x,y]$-module of rank $5$ on the elements $dx,zdx,dy,zdy,dz$. Let $w:=zdx\wedge dy-ydx\wedge dz + xdy\wedge dz$.
It follows there is an isomoprhism $\wedge^2 \Omega \cong Aw$ of $A$-modules, where $Aw$ is the free rank one $A$-module on the element $w$.
\end{lemma}
\begin{proof} The Lemma follows from an elementary calculation.
\end{proof}

One checks $\H^0_{DR}(S^2) \cong k$ is one-dimensional and also $H^2_{DR}(S^2) \cong k\overline{w}$ is one dimensional on the element $\overline{w}$. Let $\eta:=x(zdx)+y(zdy) \in \Omega$. One checks $d^1(\eta)=0$.
Let $\omega_n:=(x^2+y^2)^n\eta$ for $n \geq 0$. One also checks $d^1(\omega_n)=0$.  By definition it follows

\[  \omega_n:=(x^2+y^2)^nx(zdx) + (x^2+y^2)^ny(zdy) \in \Omega.\]

Let $a(x,y), b(x,y)z \in A:=k[x,y]\{1,z\}$ be any elements. We get the formulas

\[ d^0(a(x,y))=a_xdx+ a_ydy \]

where $a_x,a_y$ means partial derivatives with respect to the $x$ and $y$ variable.

We also get

\[  d^0(b(x,y)z)= b_x(zdx) + b_y(zdy)+ bdz \]

where $b_x,b_y$ are partial derivatives with respect to the $x$ and $y$-variables.

It follows immediately that there is no element $g(x,y,z) \in A$ with $d^0(g)=\omega_n$ for $n \geq 0$. Hence the class $\overline{\omega_n}$ is non zero in $\H^1_{DR}(S^2)$. 

I now prove that the elements $\{\overline{\omega_n}  \}_{n \geq 0}$ are linearly independent over $k$ in $\H^1_{DR}(S^2)$.

Assume $\alpha_i \in k$ and there is an equality 

\[  \sum_{i=0}^d \alpha_i \overline{\omega_i} =0 .\]

It means there is an element $g \in A$ with $d^0(g) = \sum_i \alpha_i \omega_i$. By the above calculation this implies $g=0$ and $\sum_i \alpha_i \omega_i =0$. Since the polynomials $(x^2+y^2)^n$ for $n \geq 0$ are linearly independent
over $k$ it follows $\alpha_i =0$ for all $i$ and we get the following theorem:

\begin{theorem} \label{mainalgebraic} The following holds: $dim_k(\H^i_{DR}(S^2))=1$ for $i=0,2$. The vector space $\H^1_{DR}(S^2)$ is infinite dimensional.
\end{theorem}
\begin{proof} The proof follows from the above calculation.
\end{proof}

Note: Theorem \ref{mainalgebraic} contradicts a result in \cite{grothendieck} relating algebraic deRham cohomology of a regular affine algebraic variety over the complex numbers to singular cohomology.

\section{A construction of the characteristic class}

In this section I give an explicit and elementary construction of the characteristic class $c(E) \in \Ext^1(L, \End_A(E))$. I also relate the definition to classical constructions (see \cite{teleman}, \cite{huebschmann}).

Let in the following $k$ be an arbitrary commutative unital ring and let $A$ be an arbitrary commutative unital $k$-algebra. Let $\alpha: L \rightarrow \Der_k(A)$ be an $A/k$-Lie-Rinehart algebra and let $\nabla: L \rightarrow \End_k(E)$
be an $L$-connection. Let $\nabla_1: L \rightarrow \End_k(\End_A(E))$ be defined by $\nabla_1(x)(f):=[\nabla(x), f]$ for $x\in L, f\in \End_A(E)$.  For $x,y \in L$ define the curvature $R_{\nabla}(x,y)$ as follows:

\[  R_{\nabla}(x,y):= [\nabla(x), \nabla(y)] - \nabla([x,y]) \in \End_A(E).\]

It follows we may view the curvature as an $A$-linear map  $R_{\nabla}: \wedge^2 L \rightarrow \End_A(E)$.

Let $\D(L,E, \nabla):=\End_A(E) \oplus L$ with the following $k$-bilinear product:

\[  [-,- ]: \D(L,E, \nabla) \oplus \D(L,E, \nabla) \rightarrow \D(L,E, \nabla) \]

defined by

\[  [(\phi, x), (\psi,y)] := ([\phi, \psi] + \nabla_1(x)(\psi)-\nabla_1(y)(\phi) + R_{\nabla}(x,y), [x,y])\]

for $(\phi,x), (\psi,y) \in \D(L,E,\nabla)$. There is a canonical map of $A$-modules $p_E: \D(L,E, \nabla) \rightarrow L$ defined by $p_E(\phi,x):=x$ and we get an exact sequence of $A$-modules

\[  0 \rightarrow \End_A(E) \rightarrow \D(L,E, \nabla) \rightarrow L \rightarrow 0 .\]

There is a map of $A$-modules $\alpha_E. \D(L,E,\nabla) \rightarrow \Der_k(A)$ defined by $\alpha_E(\phi,x):=\alpha(x)$ for $(\phi,x)\in \D(L,E,\nabla)$.

\begin{theorem} \label{mainconnection} Let $A$ be any commutative $k$-algebra, $(L, \alpha)$ any $A/k$-Lie-Rinehart algebra  and let $(E, \nabla)$ be any $L$-connection. The pair $(\D(L,E,\nabla), \alpha_E)$ is an $A/k$-Lie-Rinehart algebra and the map $p_E$
is a map of $A/k$-Lie-Rinehart algebras. Hence we get an exact sequence of Lie-Rinehart algebras

\[  0 \rightarrow \End_A(E) \rightarrow \D(L,E,\nabla) \rightarrow L \rightarrow 0 .\]

The map of Lie-Rinehart algebras $p_E$ has a section $s_E: L \rightarrow \D(L,E,\nabla)$  if and only if $E$ has a flat connection $\nabla_2$.

\end{theorem}  
\begin{proof} The proof will be given below.
\end{proof}

\begin{definition} \label{fundamental} The exact sequence $s(E)$ from Theorem \ref{mainconnection}

\[  0 \rightarrow \End_A(E) \rightarrow \D(L,E,\nabla) \rightarrow L \rightarrow 0 \]

is called the \emph{fundamental exact sequence} of the $L$-connection $(E, \nabla)$.

\end{definition}

\begin{example} Connections in differential geometry, algebraic geometry and algebra.\end{example}

In differential geometry a sequence similar to the sequence from Theorem \ref{mainconnection} is studied. The sequence is called the \emph{Atiyah sequence} (see \cite{huebschmann}) of a principal fiber bundle and the fundamental
exact sequence from Definition \ref{fundamental} may be seen as an algebraic version of the notion studied in differential geometry (see \cite{huebschmann}).  The history of this subject is long and complicated and I will not go into details about  this. In differential geometry  there are various types of connections such as the \emph{Ehresmann connection} and \emph{Cartan connection}. Connections in differential geometry are usually seen as splittings of relative tangent sequences of principal fiber  bundles.
These types of connections have not been studied much in the language of algebraic geometry and geometric invariant theory. 

Connections in algebra and algebraic geometry arise when studying the Riemann-Hilbert correspondence, D-modules  and p-adic cohomology. In the Riemann Hilbert correspondence one studies for a complex algebraic manifold $X$ the relation between complex representations of the topological fundamental group $\pi_1(X)$ and finite rank 
complex vector bundles $E$ on $X$ with a flat algebraic connection $\nabla$. In the theory of D-modules one studies for a regular scheme $X$ of finite type over the complex numbers, the sheaf of differential operators $D_X$ and 
the category of sheaves of modules on $D_X$. If a $D_X$-module $E$ is coherent as $\mathcal{O}_X$-module it follows there is a flat algebraic connection $\nabla$ on $E$ which "determine" $E$ as a module on $D_X$. Hence a class of $D_X$-modules
arise as pairs $(E, \nabla)$ where $\nabla$ is a flat algebraic connection on $E$ and where $E$ is a finite rank vector bundle on $X$.  Connections also show up in p-adic cohomology and in the definition of cristalline cohomology. 

Kunz has in an unpublished preprint \cite{kunz0} defined the notion "affine Riemannian scheme" and has for any affine Riemannian scheme $(X,g)$  proved there is a unique connection $\nabla$ on the tangent bundle $T_X$ compatible with the "Riemannian metric" $g$. He also develops a general theory for such "algebraic Levi-Civita connections" $\nabla$. The algebraicc Levi-Civita connection of Kunz is a Koszul-connection on the form

\[ \nabla: T_X \rightarrow \Omega^1_{X/S} \otimes T_X .\]

It is not immediate that the fundamental exact sequence is directly related to 
principal fiber bundles in algebraic geometry. One has to introduce principal fiber bundles and frame bundles using geometric invariant theory and prove the fundamental exact sequence can be constructed using the relative tangent sequence
of a principal fiber bundle. The study of principal fiber bundles and frame bundles in algebraic geometry is more involved since the construction of the quotient of a scheme by a group scheme is technically more difficult.
This problem will be studied in another paper on the subject.

Note: Under "certain conditions" there is for any $A$-module $E$ and any integer $k\geq 1$ the \emph{$k$'th order fundamental exact sequence}

\[ s_k(E):  0 \rightarrow \operatorname{S}^k(\Omega^1_{A/k})\otimes E \rightarrow J^k(E) \rightarrow J^{k-1}(E) \rightarrow 0 .\]

The sequence $s_1(E)$ is split if and only if $E$ has a \emph{Koszul connection} $\nabla_1$, which is a $k$-linear map

\[  \nabla_1: E \rightarrow \Omega^1_{A/k} \otimes_A E \]

with the property that for any $a\in A, e\in E$ it follows

\[ \nabla(ae)=a\nabla(e) + d(a) \otimes e.\]

A Koszul connection $\nabla_1$ gives rise to a $\Der_k(A)$-connection $\nabla$, which is an $A$-linear map

\[ \nabla: \Der_k(A) \rightarrow \End_k(E) \]

with the property that for any $x\in \Der_k(A), a\in A, e\in E$ it follows

\[ \nabla(x)(ae)=a\nabla(x)(e) + x(a)e.\]

Hence if $E$ is a projective module it follows $s_1(E)$ splits and hence $E$ has a Koszul connection and a corresponding $\Der_k(A)$-connection. Moreover any $\Der_k(A)$-connection $\nabla$
comes from a Koszul connection $\nabla_1$, and $R_{\nabla}=0$ if and only if $K_{\nabla}=0$. There is always a cohomology class

\[  a(E) \in \Ext^1_A(E , \Omega^1_{A/k} \otimes E)  \]

which is sometimes called the \emph{Atiyah class}. Hence there is an equality $a(E)=0$ if and only if $E$ has a Koszul connection. 

The $k$'th order fundamental exact sequence is something introduced and studied in algebraic geometry (see \cite{EGA4})
and is used in the study of ramification phenomena for morphisms of schemes and quasi coherent sheaves and the theory of discriminants. 

Given any morphism of schemes $f:X \rightarrow Y$ over a base scheme $S$, there is the relative tangent sequence

\[  0 \rightarrow T_{X/Y} \rightarrow T_{X/S} \rightarrow^{df} f^*T_{Y/S} \]

which is left exact in general. There is also the relative cotangent sequence

\[   f^* \Omega^1_{Y/S} \rightarrow \Omega^1_{X/S} \rightarrow \Omega^1_{X/Y} \rightarrow 0 \]

which is right exact in general.  In the case of the projetion map $p: P\rightarrow M$ where $P$ is a principal fiber bundle over a smooth scheme $M$, it follows the relative tangent sequence (and cotangent sequence) is exact since the projection map $p$ is smooth. The relative tangent sequence is the dual of the relative cotangent sequence and  the relative cotangent sequence is the dual of the relative cotangent sequence hence they "determine each other". They do not determine each other in general and for more general morphisms one uses the relative cotangent sequence, since this sequence has better functorial properties.

\begin{example} Left and right structures, $D$-Lie algebras and the Atiyah sequencec.\end{example}

In a recent preprint (see \cite{maa}) I have introduced the notion \emph{$D$-Lie algebra} which is a notion related to a Lie-Rinehart algebra. A $D$-Lie  algebra is essentially a Lie-Rinehart algebra $\alpha: L(f) \rightarrow \Der_k(A)$
depending on a 2-cocycle $f \in \Z^2(\Der_k(A), A)$. Using  a cocycle $f$ I construct an $A/k$-Lie-Rinehart algebra $\alpha: L(f) \rightarrow \Der_k(A)$ with the structure of a left and right $A$-module satifsying a set of conditions. This is similar to the first  order jet bundle $J^1_{A/k}(E)$ of a vector bundle $E$ - it has the structure of a left and right $A$-module. It would be interesting to study this structure in the case of the Atiyah-sequence of a principal fiber bundle
and to introduce additional structure on the classical Atiyah sequence. In \cite{maa} I study non-abelian extensions of D-Lie algebras and one wants to understand if there is a fundamental class $c^*(E)$ in this generalized setting
giving a class that is finer than the fundamental class $c(E)$ introduced in this paper.

I now prove Theorem \ref{mainconnection}.
Let us first prove that the $k$-bilinear product $[,]$ and the map $\alpha_E$ gives $\D(L,E,\nabla)$ the structure of a $A/k$-Lie-Rinehart algebra. The main "difficulty" is to prove that the product $[,]$ satisfies the Jacobi identity.

\begin{lemma} \label{derivation}   Let $\nabla_1: L \rightarrow \End_k(\End_A(E))$ be defined by $\nabla_1(x)(f):=[\nabla(x), f]$ for $x\in L, f\in \End_A(E)$. For any $x\in L, f,g\in \End_A(E)$ it follows

\[  \nabla_1(x)([f,g]) = [\nabla_1(x)(f), g ] + [f , \nabla_1(x)(g)] .\]

Hence $\nabla_1: L \rightarrow \Der_k(\End_A(E))$ and any endomorphism $\nabla_1(x)$ acts as a derivation. 
\end{lemma}
\begin{proof} The proof is an elementary calculation.
\end{proof}

Let $z_1:=(\phi, x),z_2:=(\psi, y ), z_3:=(\gamma, z) \in \D(L,E,\nabla)$. We want to prove the equality

\[   [z_1,[z_2,z_3]] +   [z_2,[z_3,z_1]]   +   [z_3,[z_1,z_2]]   =0.\]

We get the following calculation:

\[  [z_1,[z_2,z_3]]=    (  A_1, B_1 ) \in \D(L,E,\nabla) \]

where 

\[ A_1=     [\phi, [\psi, \gamma]] + [\phi, \nabla_1(y)(\gamma)] - [\phi, \nabla_1(z)(\psi) ]  + \]

\[ [\phi, R_{\nabla}(y,z) + \nabla_1(x)([\psi, \gamma]) + \nabla_1(x) \nabla_1(y)(\gamma) - \nabla_1(x)\nabla_1(z)(\psi) + \]

\[ \nabla_1(x)(R_{\nabla}(y,z)) - \nabla_1([y,z])(\phi) + R_{\nabla}(x, [y,z])\]

and 

\[ B_1= [x,[y,z]].\]

We get 

\[  [z_2,[z_3,z_1]]=    (  A_2, B_2 ) \in \D(L,E,\nabla) \]

where

\[ A_2=     [\psi, [\gamma, \phi]] + [\psi, \nabla_1(z)(\phi)] - [\psi, \nabla_1(x)(\gamma) ]  + \]

\[ [\psi, R_{\nabla}(z,x) + \nabla_1(y)([\gamma, \phi]) + \nabla_1(y) \nabla_1(z)(\phi) - \nabla_1(y)\nabla_1(x)(\gamma) + \]

\[ \nabla_1(y)(R_{\nabla}(z,x)) - \nabla_1([z,x])(\psi) + R_{\nabla}(y, [z,x])\]

and 

\[ B_2= [y,[z,x]].\]

We get

\[  [z_3,[z_1,z_2]]=    (  A_3, B_3 ) \in \D(L,E,\nabla) \]

with

\[ A_3=     [\gamma, [\phi, \psi]] + [\gamma, \nabla_1(x)(\psi)] - [\gamma, \nabla_1(y)(\phi) ]  + \]

\[ [\gamma, R_{\nabla}(x,y) + \nabla_1(z)([\phi, \psi ]) + \nabla_1(z) \nabla_1(x)(\psi) - \nabla_1(z)\nabla_1(y)(\phi) + \]

\[ \nabla_1(z)(R_{\nabla}(x,y)) - \nabla_1([x,y])(\gamma) + R_{\nabla}(z, [x,y])\]

and 

\[ B_3= [z,[x,y]].\]

We get using Lemma \ref{derivation} an equality 

\[   [z_1,[z_2,z_3]] +   [z_2,[z_3,z_1]]   +   [z_3,[z_1,z_2]]   = ( \sum A_i, \sum B_i) = 0, \]

hence the Jacobi identity holds in $\D(L,E,\nabla)$.

Let $z:=(\phi,x), w:=(\psi,y) \in \D(L,E,\nabla)$ and $a\in A$. We get

\[  [z,aw]:= [(\phi,x), (a\psi, ay)] :=\]

\[  ([\phi, a\psi] + \nabla_1(x)(\psi) - \nabla_1(ay)(\phi) + R_{\nabla}(x,ay), [x, ay]  ) =\]

\[  ( a[\phi, \psi ] + a\nabla_1(x)(\psi) + \alpha(x)(a) \psi - a\nabla_1(y)(\phi) + aR_{\nabla}(x,y), a[x,y] + \alpha(x)(a) y ) =\]

\[  a( [\phi, \psi] + \nabla_1(x)(\psi) - \nabla_1(y)(\phi) + R_{\nabla}(x,y), [x,y]) + \alpha(x)(a)(\psi, y) =\]

\[  a[z,w] + \alpha_E(z)(a) w.\]

The map $\alpha_E: \D(L,E,\nabla) \rightarrow \Der_k(A)$ is an $A$-linear map of $k$-Lie algebras hence the pair $(\D(L,E,\nabla), \alpha_E)$ is with the bracket $[,]$ an $A/k$-Lie-Rinehart algebra.
The projection map $p_E: \D(L,E,\nabla) \rightarrow L$ is an $A$-linear map of $k$-Lie algebras commuting with the maps to $\Der_k(A)$ hence the fundamental exact sequence
from Definition \ref{fundamental} is an exact sequence of Lie-Rinehart algebras. Hence I have proved Theorem \ref{mainconnection}.

Assume $s_P: L \rightarrow \D(L,E,\nabla)$ is an $A$-linear splitting of $p_E$. It follows $s_P(x):=(P(x), x)$ for a map of $A$-modules $P: L \rightarrow \End_A(E)$.   It follows $\nabla_3:=\nabla + P$ is an $L$-connection.

We get the following calculation: For any elements $x,y\in L$ it follows

\[    s_P([x,y]):=(P([x,y]), [x,y]) \]

and 

\[ [s_P(x), s_P(y)] = \]

\[  (\nabla_1(x)(P(y)) - \nabla_1(y)(P(x)) + [P(x), P(y)] + R_{\nabla}(x,y), [x,y]). \]
 
It follows

\[  \gamma_P(x,y):= [s_P(x), s_P(y)] - s_P([x,y])= (R_{\nabla_3}(x,y), 0). \]

Hence the map $s_P$ is a map of $k$-Lie algebras if and only if $\nabla_3$ is a flat connection. Hence the sequence splits in the category of $A/k$-Lie-Rinehart algebras if and only if $E$ has a flat connection.
Hence Theorem \ref{mainconnection} is proved.

The map $\gamma_P(-,-)$ gives an $A$-linear map

\[ \gamma_P: \wedge^2 L \rightarrow  \D(L,E, \nabla) \]

and if we let $P:=0$ be the zero map it follows we recover the curvature $R_{\nabla}$: We get

\[\gamma_0(x,y)=(R_{\nabla}(x,y), 0) \in \D(L,E, \nabla) .\]

Hence if we are interested in characteristic classes, it follows the study of the curvature of the $L$-connection $(E, \nabla)$ is equivalent to the study of the pair  $(\D(L,E,\nabla), \gamma_0)$. They "determine" each other.

Assume we have chosen another connection $\nabla_2:=\nabla + P$ where $P: L \rightarrow \End_A(E)$ is an $A$-linear map. We get two Lie-Rinehart algebras 

\[  \D(L, E, \nabla), \D(L, E, \nabla_1) .\]

Define the map of left $A$-modules $\phi$ as follows:

\[  \phi: \D(L,E,\nabla) \rightarrow \D(L,E,\nabla_2) \]

by

\[  \phi(f,x)=(f-P(x), x).\]

\begin{proposition} \label{iso} The map $\phi$ is an isomorphism of $A/k$-Lie-Rinehart algebras inducing an isomorphism of fundamental exact sequences: There is a commutative diagram 
\[
\diagram      0 \rto      &   \End_A(E) \rto \dto^=  &  \D(L,E,\nabla) \dto^{\phi}  \rto & L \rto \dto^= & 0    \\
                      0  \rto   &  \End_A(E)   \rto  &  \D(L,E,\nabla_2) \rto & L \rto & 0  
\enddiagram
\]

\end{proposition}
\begin{proof} The map $\phi$ is clearly an isomorphism of left $A$-modules commuting with the maps in the diagram, hence it remains to prove that the map respects the $k$-Lie product. Let $(f,x),(g,y)\in \D(L,E, \nabla)$. We get the following calculation:

\[ [ \phi(f,x), \phi(g,y)] =\]

\[  [(f-P(x), x), (g-P(y), y)] =\]

\[ ([f-P(x), g-P(y)] + \nabla_1(x)(g-P(y)) - \nabla_1(y)(f-P(x)) + R_{\nabla_2}(x,y) , [x,y]) \]

where by definition $\nabla_1(x)(\psi):=[\nabla(x), \psi]$.  A long calculation shows there is an equality

\[   [\phi(f,x), \phi(g,y)] = \]

\[  ([f,g] + \nabla_1(x)(g) - \nabla_1(y)(f) + R_{\nabla}(x,y) - P[x,y], [x,y]) =\]

\[ \phi ([(f,x), (g,y)]) \]

and it follows the map $\phi$ respects the Lie product. The Proposition follows.
\end{proof}

Let $(L, \alpha)$ be an $A/k$-Lie-Rinehart algebra where $L$ is a projective $A$-module and let $(E, \nabla)$ be an $L$-connection. Let $\Ext^1(L, \End_A(E))$ be the pointed set of equivalence classes of extensions of Lie-Rinehart algebras

\[  0 \rightarrow \End_A(E) \rightarrow L_1 \rightarrow L \rightarrow 0 .\]

Since $L$ is a projective $A$-module there is by \cite{maa} a bijection of sets

\[ \ext^1(L, \End_A(E)) \cong \Ext^1(L, \End_A(E)) \]

where $\ext^1(L, \End_A(E))$ is the following set: Let $(E, \nabla)$ be a connection and let $\nabla_1: L \rightarrow \End_k(\End_A(E))$ be the adjoint connection defined by $\nabla_1(x)(f):=[\nabla(x), f]$ for $x \in L, f \in \End_A(E)$. 
Let $S((L, \alpha), (\End_A(E), [,]) )$ be the the set of pairs $(\nabla, \psi)$ with $\nabla: L \rightarrow \Der_k(\End_A(E))$ a connection and $\psi \in \operatorname{Z}^2(L, \End_A(E))$ where for any $x,y \in L, f,g \in \End_A(E)$ it follows

\[ \nabla(x)([f,g])=[\nabla(x)(f), g] + [f , \nabla(x)(g)]  \]

and

\[ R_{\nabla}(x,y)(f)=[\psi(x,y), f] .\]

Two pairs $(\nabla, \psi), (\nabla_1, \psi_1)$ are equivalent if and only if there is a map $b \in \Hom_A(L, \End_A(E))$ with 

\begin{align}
&\label{E1}  \nabla_1(x)=\nabla(x) + [b(x), -]
\end{align}

and 

\begin{align} 
&\label{E2} \psi_1(x,y)=\psi(x,y) + d^1_{\nabla}(b)(x,y) +[b(x), b(y)].
\end{align}

\begin{lemma} Let $\nabla: L \rightarrow \Der_k(\End_A(E))$ be a connection and let $b:L \rightarrow \End_A(E))$ be an $A$-linear map. Let $\nabla_1$ be defined by $\nabla_1(x)(f):=\nabla(x)(f) + [b(x), f]$.
It follows $\nabla_1: L \rightarrow \Der_k(\End_A(E))$ is a connection with the property that for any $x \in L, f,g \in \End_A(E)$ there is an equality

\[  \nabla_1(x)([f,g])=[\nabla_1(x)(f), g] + [f , \nabla_1(x)(g)] .\]

\end{lemma}
\begin{proof} One checks the map $\nabla_1: L \rightarrow \End_k(\End_A(E))$ is an $A$-linear map. Let $x \in L$ and $f,g \in \End_A(E)$. It follows

\[  \nabla_1(x)([f,g]):=\nabla(x)([f,g]) + [b(x), [f,g]]  .\]

The Jacobi identity implies there is an equality

\[   [b(x), [f,g]]= [[b(x),f], g] + [f, [b(x),g]] .\]

It follows 

\[  \nabla_1(x)([f,g])=[\nabla_1(x)(f), g] + [f, \nabla_1(x)(g) ] .\]

 The Lemma is proved.
\end{proof}

By \cite{maa}  the equalities $\ref{E1},\ref{E2}$  define an equivalence relation on the set 

\[ S((L, \alpha), (\End_A(E), [,])) \]

 and we define $\ext^1(L, \End_A(E)):=S((L, \alpha), \End_A(E))/ \cong$ to be the quotient set. Hence we may describe the pointed set $\Ext^1(L, \End_A(E))$ in term of equivalence classes
of pairs $(\nabla, \psi)$, where $\nabla$ is a connection and $\psi$ is defined in terms of the Lie-Rinehart complex of the adjoint connection of $\End_A(E)$. This description enable us to define the structure of a torsor on $\Ext^1(L, \End_A(E))$.
One checks the adjoint connection $\nabla_1$ on $\End_A(E)$ induce a flat connection

\[ \nabla': L \rightarrow \Der_k(Z(\End_A(E)))   \]

where $Z(\End_A(E))$ is the center of the associative ring $\End_A(E)$. We may consider the Lie-Rinehart complex of the $L$-connection $(Z(\End_A(E)), \nabla')$ and the cohomology group

\[   \operatorname{H}^2(L, Z(\End_A(E))).\]

The abelian group  $ \operatorname{H}^2(L, Z(\End_A(E)))$  is the Lie-Rinehart cohomology of the flat connection  $(Z(\End_A(E), \nabla')$.

\begin{theorem} \label{mainclass} Let $(E, \nabla)$ be an $L$-connection. The fundamental exact sequence $s(E)$ defines a cohomology class $c(E) \in \Ext^1(L, \End_A(E))$ that is independent of choice of connection $\nabla$.
The class $c(E)$ is trivial if and only if $E$ has a flat connection. Let $Z(\End_A(E))$ be the center of the ring of endomorphisms $\End_A(E)$. There is an action 

\[ \sigma: \operatorname{H}^2(L, Z(\End_A(E)) \times \Ext^1(L, \End_A(E)) \rightarrow \Ext^1(L, \End_A(E)) \]

making $\Ext^1(L, \End_A(E))$ a torsor on the abelian group $\operatorname{H}^2(L, Z(\End_A(E))$.
\end{theorem}
\begin{proof}  The proof will be given in the next section.
\end{proof}

Note: Theorem \ref{mainclass} gives the pointed cohomology set $\Ext^1(L, \End_A(E))$ more structure and it may make it easier to calculate $\Ext^1(L, \End_A(E))$  "by hand" using the techniques introduced in this paper.
If $E$ is a finite rank locally trivial module, it follows $Z(\End_A(E)) \cong A$ and hence $\operatorname{H}^2(L, Z(\End_A(E)) \cong \operatorname{H}^2(L,A)$ and the latter is easier to calculate.

\begin{definition} The cohomology class $c(E) \in\Ext^1(L, \End_A(E))$ defined in Theorem \ref{mainclass} is the \emph{fundamental class} of the $A$-module $E$.
\end{definition}

Note: By Proposition \ref{iso} and Theorem \ref{mainclass} it follows the class $c(E)$ is independent of choice of connection. The fundamental class $c(E)$ lives in a cohomology torsor and the Atiyah class $a(E)$ lives in a cohomology group,
hence these two classes differ. The class $a(E)=0$ if and only if $E$ has a connection. If $a(E)=0$ it follows $c(E)$ is trivial if and only if $E$ has a flat connection.


\begin{example} The classical construction of Teleman and the Chern character. \end{example}

In the paper \cite{teleman}, Teleman constructs for any extension of $A/k$-Lie-Rinehart algebras

\[  0 \rightarrow W \rightarrow L_1 \rightarrow L \rightarrow 0 \]

an any "invariant function"  $g \in \operatorname{I}^i_A(W, A)$ a characteristic class 

\[ [g(\Omega_{\gamma})] \in \H^{2i}(L,A) \] 

that is independent of "choice of connection". Here $\H^{2i}(L,A)$ is the Lie-Rinehart cohomology group of $(L, \alpha)$. 
A "connection" in the language of Teleman is an $A$-linear splitting $\gamma: L \rightarrow L_1$ of the exact sequence. 
Hence the characteristic class of Teleman lives in a cohomology group. The characteristic class of Teleman may be seen as a characteristic class of an arbitrary "principal fiber bundle" $\pi: P_1 \rightarrow S$, where we view the map  $p: L_1 \rightarrow L$  as the relative tangent sequence of the projection map $\pi$.  The class introduced in this paper lives in a pointed cohomology torsor, hence these classes differ.

If $W:=\End_A(E)$ with $(E, \nabla)$ a connection,  it follows the characteristic ring introduced in \cite{teleman}
is trivial if $E$ has a flat connection. This is proved in Corollary 12 in \cite{teleman}. Hence if the class $c(E)$ from this section is trivial it follows the class $[g(\Omega)]$ and the characteristic ring in \cite{teleman} is trivial. 
Hence in the case when $W:=\End_A(E)$ it follows the class $c(E)$ can be used to study the classes and the characteristic ring of Teleman since the curvature of the connection $\nabla$ on $E$ is completely determined by the splitting $\gamma$.

Note: Using  the fundamental exact sequence

\[ 0 \rightarrow \End_A(E) \rightarrow \D(L, E, \nabla) \rightarrow L \rightarrow 0 \]

and any invariant function     $g \in \operatorname{I}^i_A(\End_A(E), A)$,  we get using Teleman's construction a cohomology class $[g(\Omega(E))] \in \H^{2i}(L,A)$ that is independent  of choice of connection $\nabla$.
If the class $[g(\Omega(E))]$ satisfies the Whitney sum formula it will be trivial for stably trivial vector bundles.  When $g$ is the trace map we get the Chern classes $c_i(E) \in \H^{2i}(L,E)$ as introduced in \cite{maa1}.
The Chern class $c_i(E)$ from \cite{maa1} is zero for stably trivial vector bundles since it satisfies the Whitney sum formula.

\section{The torsor structure on $\Ext^1(L, \End_A(E))$.}

In this section I define and prove some properties of the torsor structure on the pointed cohomology set $\Ext^1(L, \End_A(E))$. 

Let $k$ be a commutative unital ring and let $A$ be a commutative unital $k$-algebra and let $(L, \alpha)$
be an $A/k$-Lie-Rinehart algebra. Let $\nabla: L \rightarrow \End_k(E)$ be an $L$-connection and let $\nabla_1: L \rightarrow \End_k(\End_A(E))$ be the adjoint connection defined as follows: Let $x\in L, f \in \End_A(E)$ and define 
$\nabla_1(x)(f):=[\nabla(x), f]$. The adjoint connection $\nabla_1$ has the following property: For any $x,\in L, f,g \in \End_A(E)$ it follows

\[ \nabla_1(x)([f,g])=[\nabla_1(x)(f), g] + [f, \nabla_1(x)(g)]. \]

If $W$ is an $A$-module with the structure of an $A$-Lie algebra with $A$-linear  bracket $[,]$ satisfying the Jacobi identity and if $\nabla_2: L \rightarrow \End_k(W)$ is an $L$-connection with the property $\ref{star}$:  For any $x \in L, u,v \in W$ there is an equality

\begin{align}
\label{star}   \nabla_2(x)([u,v]) =[\nabla_2(x)(u), v] + [u, \nabla_2(x)(v)].
\end{align}

We write this as follows: The $A$-linear map $\nabla_2$ is a map of $A$-modules

\[ \nabla_2: L \rightarrow \Der_k(W) .\]

\begin{definition} Let $(W, [,])$ be an $A$-Lie algebra and let $\nabla: L \rightarrow \End_k(W)$ be an $L$ connection satifsying $\ref{star}$ above. We say $\nabla$ is a \emph{connection on the $A$-Lie algebra $W$} and we write

\[\nabla: L \rightarrow \Der_k(W)\]

to indicate this.

\end{definition}

If $\nabla$ is an $L$-connection on the $A$-Lie algebra $(W, [,])$ it follows the curvature $R_{\nabla}$ of $\nabla$  has the following property: For any $x,y\in L, w,v \in W$ it follows

\[  R_{\nabla}(x,y)([v,w])= [  R_{\nabla}(x,y)(v), w] + [v,    R_{\nabla}(x,y)(w)].\]

The "module of derivations" $\Der_k(W) \subseteq \End_k(W)$ is the sub-$A$-module of $k$-linear endomorphism $\rho \in \End_k(W)$ satisfying the following "derivation property": For any $u,v \in W$ the following holds:

\[ \rho([u,v])=[\rho(u), v] + [u, \rho(v)] .\]

The module $\End_A(E)$ has a bracket $[,]$ making it an $A$-Lie algebra, and it follows the adjoint connection $\nabla_1$ gives an $A$-linear map

\[ \nabla_1: L \rightarrow \Der_k(\End_A(E)) .\]

Consider the construction from  \cite{maa}.   Let $(W, [,])$ be an arbitrary $A$-Lie algebra and let $S((L, \alpha), (W, [,])]$ be the set of pairs $(\nabla, \psi)$ with $\nabla: L \rightarrow \Der_k(W)$ a connection and $\psi \in \Z^2(L,W)$
a 2-cocycle with the property that for any $x,y\in L, w\in W$ it follows there is an equality

\[  R_{\nabla}(x,y)(w)=[\psi(x,y), w].\]

By definition $R_{\nabla}(x,y):=[\nabla(x), \nabla(y)]- \nabla([x,y])$.  In \cite{maa} I define an equivalence relation on $S:=S((L, \alpha), (W, [,]))$ as follows: We say two elements $(\nabla_1, \psi_1), (\nabla_2, \psi_2)$ in $ S$
are equivalent if and only if there is an element $b\in \Hom_A(L,W)$ with the following property: Let  $P_b \in \Hom_A(L, \End_A(W))$ be the follwing element: For any $x\in L, w\in W$ define

\[  P_b(x)(w):=[b(x), w] .\]

We say $(\nabla_1 , \psi_1) \cong (\nabla_2, \psi_2)$ if and only if there are equalities

\[  \nabla_2(x)=\nabla_1(x) + P_b(x) \]

and

\[  \psi_2(x,y)=\psi_1(x,y) + d^1_{\nabla_1}(b)(x,y) + [b(x), b(y)].\]

\begin{lemma} The map $P: \Hom_A(L,W) \rightarrow \Hom_A(L, \End_A(W))$ defined by $P$ has the following properties: Let $a,b \in \Hom_A(L,W)$ and let $u,v \in W$. There are equalities $P_{a+b}=P_a + P_b$ and $P_a([u,v])=[P_a(u), v] +[u, P_a(v)]$.
Let $\nabla: L \rightarrow \Der_k(W)$ be a connection satisfying the property $\ref{star}$ and let $b \in \Hom_A(L,W)$. It follows $\nabla_1:=\nabla + P_b$ is a connection satisfying property $\ref{star}$.
\end{lemma}
\begin{proof} The proof is a straight forward calculation.
\end{proof}

I will now for the sake of "completeness" prove that the relation defined above is an equivalence relation. The relation is clearly reflexive and I will prove it is symmetric and transitive. 

Assume $(\nabla_1, \psi_1) \cong (\nabla_2, \psi_2)$. By definition there is an element $a\in \Hom_A(L,W)$ with the property that

\[  \nabla_2(x)=\nabla_1(x) + P_a(x) \]

and

\[  \psi_2(x,y)=\psi_1(x,y) + d^1_{\nabla_1}(a)(x,y) + [a(x), a(y)].\]

Let $b:=-a \in \Hom_A(L,W)$. It follows immediatley there is an equality

\[  \nabla_1(x)=\nabla_2(x)-P_a(x)=\nabla_2(x)+P_b(x) .\]

By assumtion there is an equality

\[   \psi_1(x,y) = \psi_2(x,y) - d^1_{\nabla_1}(a)(x,y) -[a(x),a(y)].\]

We get the following:

\[ \psi_2(x,y)+ d^1_{\nabla_2}(b)(x,y)+[b(x),b(y)] =\]

\[  \psi_2(x,y) + \nabla_2(x)(b(y))-\nabla_2(y)(b(x)) + a([x,y])  +[a(x),a(y)]   =\]

\[  \psi_2(x,y) - \nabla_1(x)(a(y)) - [a(x), a(y)] +\nabla_1(y)(a(x)) +[a(y), a(x)] + a([x,y]) + [a(x), a(y)] =\]

\[  \psi_2(x,y) - \nabla_1(x)(a(y)) + \nabla_1(y)(a(x)) + a([x,y]) + [a(y), a(x)] =\]

\[ \psi_2(x,y) - d^1_{\nabla_1}(a)(x,y) -[a(x), a(y)]  =\psi_1(x,y).\]

hence it follows $(\nabla_2, \psi_2) \cong (\nabla_1, \psi_1)$ and the relation is symmetric. I now prove the relation is transitive.   Assume $(\nabla_i, \psi_i) \in S$ for $i=1,2,3$ are equivalent elements in $S$. There are elements $a,b \in \Hom_A(L,W)$ 
and equalities

\[  \nabla_2(x)=\nabla_1(x) + P_a(x) \]

and

\[  \psi_2(x,y)=\psi_1(x,y) + d^1_{\nabla_1}(a)(x,y) + [a(x), a(y)],\]

and

\[  \nabla_3(x)=\nabla_2(x) + P_b(x) \]

and

\[  \psi_3(x,y)=\psi_2(x,y) + d^1_{\nabla_2}(b)(x,y) + [b(x), b(y)].\]

Let $c:=a+b$ and consider the element $P_c$. By definition we get equalities

\[  \nabla_3 = \nabla_2 + P_b = \nabla_1 + P_a + P_b = \nabla_1 + P_c.\]

We get

\[  \psi_3(x,y)= \psi_2(x,y) + d^1_{\nabla_2}(b)(x,y) + [b(x), b(y)] =\]

\[  \psi_3(x,y)= \psi_1(x,y) + d^1_{\nabla_1}(a)(x,y)+ [a(x), a(y)  ] + d^1_{\nabla_2}(b)(x,y) + [b(x), b(y)] =\]

\[ \psi_1(x,y) + \nabla_1(x)(a(y)) - \nabla_1(y)(a(x)) -a([x,y]) + [a(x),a(y)] + \]

\[ \nabla_2(x)(b(y)) - \nabla_2(y)(b(x)) -b([x,y]) + [b(x), b(y)] =\]

\[ \psi_1(x,y) +  \nabla_1(x)(c(y) - \nabla_1(y)(c(x)) - c([x,y]) + \]

\[  [a(x), a(y)] + [a(x), b(y)] + [b(x), a(y)] + [b(x), b(y)] =\]

\[ \psi_1(x,y) + \nabla_1(x)(c(y)) - \nabla_1(y)(c(x)) - c([x,y]) + [c(x), c(y)] .\]

\[ \psi_3(x,y)= \psi_1(x,y) + d^1_{\nabla_1}(c)(x,y) + [c(x), c(y)] .\]

It follows $(\nabla_1, \psi_1) \cong (\nabla_3, \psi_3)$ and the relation is transitive.  

\begin{definition}  Let  $\ext^1(L, \End_A(E)):=S((L, \alpha), (\End_A(E),[,]))/ \cong$ be the quotient.
\end{definition}

In Corollary 2.13 in \cite{maa} I prove there is a one to one correspondence of sets

\[  \Ext^1(L, \End_A(E)) \cong \ext^1(L, \End_A(E)) \]

when $L$ is projective as left $A$-module. I will use this description of the pointed cohomology set $\Ext^1(L, \End_A(E))$ to define the structure of a torsor on the set $\Ext^1(L, \End_A(E))$ of equivalence classes of extensions of 
Lie-Rinehart algebras.

Let 

\[  0 \rightarrow W \rightarrow L_1 \rightarrow^p L \rightarrow 0 \]

be an exact sequence of $A/k$-Lie-Rinehart algebras $(L_1, \alpha_1), (L, \alpha)$ and let $s: L \rightarrow L_1$ be a fixed $A$-linear section of $p$. Define the map $\nabla_s \in \Hom_A(L, \Der_k(W))$ as follows:

\[  \nabla_s(x)(w):=[s(x), w] \]

where $[,]$ is the bracket on $L_1$ and where $x\in L, w\in W$. Define the map $\psi_s: \wedge^2 L \rightarrow  W$ by

\[ \psi_s(x,y):=[s(x), s(y)]-s([x,y]) \]

where $x,y \in K$.  It follows $\psi_s \in \Z^2(L,W)$ where we use the connection $\nabla_s$ to define the differential. The section $s$ defines a splitting $L_1 \cong W \oplus L$ with induced map

\[  p: W \oplus L \rightarrow L \]

defined by $p(w,x):=x$ and $s(x):=(0,x) \in W \oplus L$ for $(w,x) \in L_1, x \in L$. Using this decomposition of $L_1$ into a direct sum it follows the bracket $[,]$ looks as follows:
For any $(w,x), (v,y)$ we get 

\[  [ (w,x), (v,y)]=([w,v] + \nabla_s(x)(v) - \nabla_s(y)(w) + \psi_s(x,y), [x,y]_L ) \]

where $[,]_L$ is the Lie bracket on $L$.  Let $T((L, \alpha),s, (W,[,]))$ be the set of pairs $(\nabla, \psi)$ with 

\[ \nabla: L \rightarrow \Der_k(W) \]

an $L$-connection   and $\psi: \Z^2(L,W)$  a 2-cocycle with the property there is an element $b \in \Hom_A(L,W)$ such that there are equalities

\[  \nabla(x) = \nabla_s(x) + P_b(x) \]

and 

\[  \psi(x,y)= \psi_s(x,y) + d^1_{\nabla_s}(b)(x,y) + [b(x), b(y)] .\]

\begin{proposition} \label{bijection} There is a one-to-one correspondence between the set of $A$-linear sections of the map $p:L_1 \rightarrow L$ and the set $T((L, \alpha), s,(W, [,]))$.
\end{proposition}
\begin{proof} Given an $A$-linear  section $t$ of $p$, in the decomposition $L_1 \cong W \oplus L$ we must have $t(x)=(b(x), x)$ for some $b \in \Hom_A(L,W)$ and one checks there are equalities

\[  \nabla_t(x) = \nabla_s(x) + P_b(x) \]

and 

\[  \psi_t(x,y)= \psi_s(x,y) + d^1_{\nabla_s}(b)(x,y) + [b(x), b(y)] .\]

Conversely given a pair $(\nabla, \psi)$ in $T((L, \alpha), s, (W, [,]))$and an element $b \in \Hom_A(L,W)$  and equalities 

\[  \nabla(x) = \nabla_s(x) + P_b(x) \]

and 

\[  \psi(x,y)= \psi_s(x,y) + d^1_{\nabla_s}(b)(x,y) + [b(x), b(y)] .\]

Define the section $t(x):=(b(x), x)$ of $p$.  It follows $\nabla_t = \nabla, \psi_t=\psi$. Hence the two sets are in bijection and the Proposition is proved.
\end{proof}

\begin{example} Characteristic classes for principal fiber bundles and the Teleman construction. \end{example}

Fix a section $s$ of the projection map $p: L_1 \rightarrow L$.
Teleman constructs for any section $\gamma$ of the projection map $p$ a characteristic class. The section $\gamma$ corresponds by Proposition \ref{bijection} one to one to a pair $(\nabla_{\gamma}, \psi_{\gamma})$ of maps 

\[ \nabla_{\gamma}: L \rightarrow \Der_k(W) \]

and 

\[ \psi_{\gamma}: \wedge^2 L \rightarrow W \]

where $\nabla_{\gamma}$ is a connection, $\psi_{\gamma} \in \Z^2(L, W)$ and where for any $x,y \in L, w\in W$ it follows 

\[  R_{\nabla_{\gamma}}(x,y)(w)= [\psi_{\gamma}(x,y), w].\]

Moreover there is an element $b\ in \Hom_A(L,W)$ with 

\[  \nabla_{\gamma}(x) = \nabla_s(x) + P_b(x) \]

and 

\[  \psi_{\gamma}(x,y)= \psi_s(x,y) + d^1_{\nabla_s}(b)(x,y) + [b(x), b(y)] .\]

Hence Teleman's characteristic class introduced in \cite{teleman} is a characteristic class of the triple $(W, \nabla_{\gamma}, \psi_{\gamma})$, or equivalently the connection $(W, \nabla_{\gamma})$ since the map $\psi_{\gamma}$ is uniquely 
determined by the curvature $R_{\nabla_{\gamma}}$ of $\nabla_{\gamma}$. Since the class introduced by Teleman is independent of choice of connection $\gamma$, it is a characteristic class $[g(W)]$ depending only on the 
polynomial $g$ and the $A$-Lie algebra  $(W,[,])$. 

Associated to the connection $(W, \nabla_{\gamma})$ there is a fundamental exact sequence

\[  0 \rightarrow \End_A(W) \rightarrow \D(L, W, \nabla_{\gamma}) \rightarrow L \rightarrow 0 \]

and a cohomology class $c(W) \in \Ext^1(L, \End_A(W))$, and $c(W)$ is trivial if and only if $W$ has a flat connection. By Corollary 16 in \cite{teleman}: If $W$ has a flat connection it follows $[g(W)]$ is trivial.
Hence the fundamental class $c(W)$ of $W$ can be used in the study of Teleman's class in complete generality. 

A large number of people have contributed to the  study of Atiyah sequences and Lie algebroids in differential geometry and Lie-Rinehart algebras in commutative algebra and some names to be mentioned are Coste, Dazord, Pradines and Weinstein (see \cite{mackenzie2}).  Mackenzie, Huebschmann and Kubarski have constructed similar exact sequences for  "principal fiber bundles" in an algebraic and differential geometric framework in \cite{huebschmann0} and \cite{kubarski} in some special cases. In \cite{huebschmann} and \cite{kubarski} and \cite{mackenzie}
you find many references and also some history on the topic.

\begin{example} The Euler class of an inner product module. \end{example}

In \cite{kong} Kong defines the \emph{Euler class} $e(P, h)$ of an inner product module $(P,h)$, and he proves that the class $e(P,h)$ detects the non triviality of the tangent bundle on the real $2n$-sphere. Kong proves in
\cite{kong} that for any  finite rank projective $A$-module $P$ with an inner product $h$, there is a unique connection $\nabla$ compatible with $h$. Then he uses the Pfaffian $Pf$ and the curvature $K_{\nabla}$ of $\nabla$
to construct $e(P,h)$. He also proves that the class $e(P,h)$ is independent of choice of connection. If $P$ has a flat connection $\nabla$ compatible with $h$,  it follows $e(P,h)$ is trivial.
Since $P$ is a finite rank and projective $A$-module we may construct the class $c(P)$ and if $c(P)$ is trivial it follows $e(P,h)$ i trivial, hence the fundamental classs $c(P)$ can be used to study $e(P,h)$.

\begin{theorem} Let $\nabla: L \rightarrow \End_k(E)$ be any connection where $L$ is a projective $A$-module. Let $\nabla_1: L \rightarrow \Der_k(\End_A(E))$ be the adjoint connection. Let $Z(\End_A(E))$ be the center of the associative ring $\End_A(E)$. The adjoint connection $\nabla_1$ induce a flat connection $\nabla_2: L \rightarrow \Der_k(Z(\End_A(E)))$. There is an action of the abelian group $\H^2(L, Z(\End_A(E)))$

\[ \sigma :\H^2(L, Z(\End_A(E))) \times \Ext^1(L, \End_A(E)) \rightarrow \Ext^1(L, \End_A(E)) \]

giving $\Ext^1(L, \End_A(E))$ the structure of a torsor on $\H^2(L, Z(\End_A(E)))$.
\end{theorem}
\begin{proof}  One checks there is an induced connection $\nabla_2$ on $Z(\End_A(E))$ that is flat.  Since $L$ is projective it follows  there is by \cite{maa} a bijection of sets

\[ \Ext^1(L, \End_A(E)) \cong \ext^1(L, \End_A(E)) .\]

Choose any representative $(\nabla, \psi)$ of a class in $\ext^1(L, \End_A(E))$ and define for any representative $\rho \in \Hom_A(\wedge^2, Z(\End_A(E)))$ the following: Let $\sigma(\rho, (\nabla, \psi))$ be the equivalence class of the element $(\nabla, \psi+ \rho)$. One easily checks this gives a well defined  action $\sigma$ and the Theorem follows.
\end{proof}

Note: When $E$ is a finite rank projective $A$-module it follows there is an equality $Z(\End_A(E)) \cong A$ hence $\H^2(L, Z(\End_A(E))) \cong \H^2(L,A)$, and the latter is a non-trivial abelian group in general.
The torsor structure on $\Ext^1(L, \End_A(E))$ may help when calculating $\Ext^1(L, \End_A(E))$ in explicit examples.

If $k$ is a field of characteristic $p>0$ and $k \subseteq L$ is a purely inseparable extension of exponent one, it follows restricted Lie-Rinehart cohomology $\H^2_{res}(\Der_k(L), L)$ of the restricted Lie-Rinehart algebra $\Der_k(L)$
calculates the Brauer group $\operatorname{Br}^L_k$ (see \cite{dokas} and \cite{hochschild}).  Hence if there is a "restricted version" of the pointed cohomology set $\Ext^1(L,\End_A(E))$ in characteristic $p>0$,
it may be this set is a torsor on a  Brauer group. There is a large litterature on connections in characteristic $p>0$ and its relation to algebraic K-theory and "motives" (see \cite{gille}).


\begin{thebibliography}{4}


\bibitem{bourbaki1} N. Bourbaki,  Éléments de mathématique. Algèbre. Chapitres 1 à 3, Springer Verlag (2007)

\bibitem{bourbaki2} N. Bourbaki,  Elements of mathematics. Algebra II. Chapters 4–7 Springer Verlag (2003)

\bibitem{DG} M. Demazure, Michel, P. Gabriel, Algebraic groups. Volume I: Paris: Masson et Cie, Éditeur; Amsterdam: North-Holland Publishing Company. xxvi, 700 p. (1970). 

\bibitem{gille} Gille, P, Szamuely, T, Central simple algebras and Galois cohomology. 2nd ed, Cambridge Studies in Advanced Mathematics 165. Cambridge: Cambridge University Press 417 p. (2017). 


\bibitem{EGA4} A. Grothendieck, Elements de geometrie algebrique IV4, Etude locale des schemas et des morphismes de schemas , Publ. Math. IHES no. 32 (1967)

\bibitem{SGA3I} M. Demazure, A. Grothendieck, 
Schémas en groupes. I: Propriétés générales des schémas en groupes. Exposés I à VIIb. Séminaire de Géométrie Algébrique 1962/64,  Lecture Notes in Mathematics. 151. Springer-Verlag  (1970). 

\bibitem{SGA3II} M. Demazure, A.  Grothendieck, Schémas en groupes. II: Groupes de type multiplicatif, et structure des schémas en groupes généraux. Exposés VIII à XVIII. Séminaire de Géométrie Algébrique du Bois Marie 1962/64 (SGA 3) 
Lecture Notes in Mathematics. 152. Springer-Verlag  (1970). 

\bibitem{SGA3III} M. Demazure, A.  Grothendieck,  Schémas en groupes. III: Structure des schémas en groupes réductifs. Exposés XIX à XXVI. Séminaire de Géométrie Algébrique du Bois Marie 1962/64 (SGA 3),  
Lecture Notes in Mathematics. 153.  Springer-Verlag (1970). 


\bibitem{docarmo}Do Carmo, M,  Riemannian geometry, Mathematics: Theory and Applications,  Birkhäuser  (1992). 


\bibitem{dokas} Dokas, I, Cohomology of restricted Lie-Rinehart algebras and the Brauer group, Adv. Math. 231, No. 5, 2573-2592 (2012). 

\bibitem{grothendieck} Grothendieck, A, On the De Rham cohomology of algebraic varieties, Publ. Math., Inst. Hautes Étud. Sci. 29, 95-103 (1966). 

\bibitem{hochschild} Hochschild, G, Simple algebras with purely inseparable splitting fields of exponent 1, Trans. Am. Math. Soc. 79, 477-489 (1955). 

\bibitem{huebschmann0} Huebschmann, J, Extensions of Lie-Rinehart algebras and the Chern-Weil construction, Higher homotopy structures in topology and mathematical physics, Contemp. Math. 227, 145-176 (1999). 


\bibitem{huebschmann} Huebschmann, J. Poisson cohomology and quantization,  \emph{arXiv:1303.3903}


\bibitem{kobayashi} Kobayashi, S. , Katsumi N. Foundations of differential geometry. I. Wiley  (1963).



\bibitem{kong} Kong, M, Euler classes of inner product modules, J. Algebra 49, 276-303 (1977). 


\bibitem{kraft} H. Kraft, Challenging problems on affine n-space, Sem. Bourbaki exp. no 802, Asterisque tome 237 (1996) 


\bibitem{kubarski} Kubarski, J, Characteristic classes of regular Lie algebroids – a sketch, The proceedings of the 11th winter school on geometry and physics held in Srní, Czechoslovakia, Suppl. Rend. Circ. Mat. Palermo, II. Ser. 30, 71-94 (1993). 


\bibitem{kunz0} Kunz, E, Algebraic differential calculus, \emph{unpublished preprint} (2002)

\bibitem{kunz} Kunz, E, Kähler differentials, Advanced Lectures in Mathematics   (1986). 

\bibitem{loday} Loday, J.L.  Cyclic homology. 2nd ed., Grundlehren der Mathematischen Wissenschaften. 301. Berlin: Springer. xviii, 513 p. (1998). 

\bibitem{maa1} Maakestad, H,  The Chern character for Lie-Rinehart algebras, Ann. Inst. Fourier 55, No. 7, 2551-2574 (2005). 


\bibitem{maa} Maakestad, H, Extensions of Lie algebras of differential operators \emph{arXiv:1512.02967} (2019)


\bibitem{mackenzie} Mackenzie, K, Lie groupoids and Lie algebroids in differential geometry, London Mathematical Society Lecture Note Series 124 (1987). 

\bibitem{mackenzie2} Mackenzie, K, The general theory of Lie groupoids and Lie algebroids, London Mathematical Society Lecture Note Series 213  (2005). 


\bibitem{matsumura} Matsumura, H,  Commutative ring theory, Cambridge Studies in Advanced Mathematics, 8  Cambridge University Press  (1989). 

\bibitem{teleman} Teleman, N, A characteristic ring of a Lie algebra extension. II, Atti Accad. Naz. Lincei, VIII. Ser., Rend., Cl. Sci. Fis. Mat. Nat. 52, 708-711 (1972). 


\end{thebibliography}
\end{document}